\documentclass[12pt,a4paper]{article}
\usepackage{amsmath,amssymb,amsthm}
\usepackage[dvips,final]{graphicx}

\usepackage{color}
\usepackage{xcolor} 
\usepackage{soul}

\definecolor{vio}{rgb}{0.5,0,0.4}
\definecolor{gre}{rgb}{0.1,0.6,0}
\definecolor{gray}{rgb}{0.82,0.82,0.82}
\definecolor{grayl}{rgb}{0.9,0.9,0.9}

\newenvironment{pkre}{\color{blue}}{\color{black}}
\newcommand{\bpk}{\begin{pkre}}
\newcommand{\epk}{\end{pkre}}

\newenvironment{pkred}{\color{magenta}}{\color{black}}
\newcommand{\br}{\begin{pkred}}
\newcommand{\er}{\end{pkred}}

\newtheorem{theorem}{Theorem}[section]

\newtheorem{lemma}[theorem]{Lemma}
\newtheorem{corollary}[theorem]{Corollary}
\newtheorem{proposition}[theorem]{Proposition}
\newtheorem{problem}[theorem]{Problem} 
\theoremstyle{definition}
\newtheorem{definition}[theorem]{Definition}
\newtheorem{remark}[theorem]{Remark}

\numberwithin{equation}{section}

\usepackage{titlesec}
\titleformat{\section}{\bfseries}{\thesection}{1em}{}
\titleformat{\subsection}{\itshape}{\thesubsection}{1em}{}

\newcommand{\real}{\mathbb{R}}
\newcommand{\nat}{\mathbb{N}}

\newcommand{\dis}{\displaystyle}
\newcommand{\Var}{\mathop{\mathrm{Var}}}
\def\expe{\mathrm{e}}
\def\dd{\,\mathrm{d}}
\def\dist{\mathrm{dist}}

\def\ok{^{(k)}}
\def\ol{^{(l)}}

\def\scal#1{\left\langle #1 \right\rangle}
\def\ve{\varepsilon}

\def\for{\mbox{ for }\ }

\def\be{\begin{equation}\label}
\def\ee{\end{equation}}

\newfont{\ctv}{msam10}

\newcommand{\bbox}{\mbox{\ctv \symbol{4}}}
\def\QED{{${}\hfill\bbox$}}
\newenvironment{pf}[1]{\par\vskip1mm{\noindent\it #1.}\ }{\QED\par
\vskip2mm}

\def\bpf{\begin{pf}}
\def\epf{\end{pf}}

\begin{document}

\title{Non-convex sweeping processes in the space of regulated functions\footnote{Supported by the GA\v CR Grant No.~20-14736S, RVO: 67985840, and by the European Regional Development Fund, Project No. CZ.02.1.01/0.0/0.0/16{\_}019/0000778.}}

\date{}

\author{Pavel Krej\v c\'\i\footnote{Faculty of Civil Engineering, Czech Technical University, Th\'akurova 7, 16629 Praha 6, Czech Republic (\texttt{pavel.krejci@cvut.cz}).}, $\;$ 
Giselle Antunes Monteiro\footnote{Institute of Mathematics, Czech Academy of Sciences, \v Zitn\'a 25,
11567~Praha 1, Czech Republic,  (\texttt{gam@math.cas.cz}).}, $\;$ 
Vincenzo Recupero\footnote{Dipartimento di Scienze Matematiche, Politecnico di Torino, Corso Duca degli Abruzzi 24, 10129 Torino, Italy (\texttt{vincenzo.recupero@polito.it}).} \footnote{Vincenzo Recupero is a member of GNAMPA-INdAM.}}


\maketitle

\begin{abstract}
The aim of this paper is to study a wide class of non-convex sweeping processes with moving constraint whose translation and deformation are represented by regulated functions, i.\,e., functions of not necessarily bounded variation admitting one-sided limits at every point. Assuming that the time-dependent constraint is uniformly prox-regular and has uniformly non-empty interior, we prove existence and uniqueness of solutions, as well as continuous data dependence with respect to the sup-norm.

\ 

\begin{footnotesize}
\end{footnotesize}
\end{abstract}

\section*{Introduction}

Sweeping processes were introduced in \cite{Mor71} as an abstract setting of problems arising for example in elastoplasticity modeling, where the constitutive relation can be formulated as a constrained evolution system. Typically, the functional framework consists in assuming that
\be{i1}
X \ \text{ is a real Hilbert space }
\ee
endowed with scalar product $\scal{\cdot, \cdot}$ and norm $|x| = \sqrt{\scal{x,x}}$ for $x \in X$, and one considers a family $C(t) \subset X$ of nonempty closed subsets of $X$ parameterized by the time variable $t \in [0,T]$, where $T>0$ is some given final time. The problem is to find a function 
$\xi: [0,T] \to X$ with a prescribed initial condition $\xi(0) = \xi_0 \in C(0)$, such that $\xi(t) \in C(t)$ for all 
$t \in [0,T]$ and its derivative at time $t$ points in the outward normal direction to $C(t)$ at the point 
$\xi(t)$. Formally, this can be stated as
\be{i1a}
-\dot \xi(t) \in N_{C(t)}(\xi(t)) \quad \for t \in (0,T), \quad \xi(0) = \xi_0,
\ee
where both the ``time derivative" $\dot \xi(t)$ and the outward normal cone $N_{C(t)}(\xi(t))$ to $C(t)$ at the point 
$\xi(t)$ have to be given an appropriate meaning.

In the paper \cite{Mor71}, this problem is uniquely solved provided that $C(t)$ is convex for every time $t$ and that the mapping $t \mapsto C(t)$ is absolutely continuous in terms of the Hausdorff distance. In this case the solution 
$\xi$ turns out to be absolutely continuous and \eqref{i1a} is satisfied almost everywhere. In \cite{moreau} the analysis of sweeping processes was then extended to the case when the convex moving set $C(t)$ has bounded variation with respect to the Hausdorff metric. Under this weaker assumption, inclusion \eqref{i1a} has to be properly interpreted in the sense of the differential measures and it is shown to admit a unique solution of bounded variation.

The technique introduced in \cite{moreau} is based on the so-called \emph{catching up algorithm} and consists in approximating $C(t)$ by a sequence of right continuous step convex-valued functions 
$C_k(t)$, i.\,e., functions such that $[0,T]$ is partitioned into a finite number of intervals where $C_k(t)$ is constant. The approximate solution is constructed as a step function $\xi_k(t)$ by an iterative process, where the next value is obtained by projection onto the current set $C_k(t)$. The argument then consists in proving that the sequence $\{\xi_k\}$ uniformly converges to a right continuous $BV$ function $\xi$ solving the suitable generalized version of \eqref{i1a}, that can be also represented by the integral variational inequality (see also \cite{RecSan18})
\begin{equation}\label{i5}
  \int_0^T \scal{\xi(t) - z(t), \dd \xi(t)} \le 0 \qquad \text{for all $z : [0,T] \to X$, $z(t) \in C(t)$}, 
\end{equation}
where the test functions $z$ are required to have some regularity properties, e.\,g., bounded variation, and where the integral is understood in terms of the differential measure $\dd \xi$. 

A relevant particular case of sweeping process occurs when the constraint $C(t)$ has a fixed shape and moves only by means of translation, i.\,e., if $C(t)$ is of the form $C(t) = u(t) - Z$ for a given function 
$u : [0,T] \to X$ and a fixed closed convex set $Z \subset X$. The resulting input-output relation 
$u \mapsto \xi$ is called the {\em (vector) play operator} and it was independently studied in the monograph \cite{KraPok} when $X$ is finite dimensional, $Z$ is bounded with non-empty interior, and $u$ is continuous. An extension to the space of regulated functions has been done in \cite{KreLau02} and \eqref{i5} is understood in the sense of Kurzweil or Young integral. Note that the Young integral can be interpreted as a variant of the Kurzweil integral, see \cite{negli}. A comparison between the measure approach and the Kurzweil/Young integral approach to \eqref{i5} is discussed in 
\cite[Section A.4]{Rec11a}.

Indeed, the integral in \eqref{i5} makes clear sense only if the solution $\xi$ is of bounded variation. This can be achieved if the moving constraint has non-empty interior, as it has been shown for the case of continuous inputs independently in \cite{Cas83} and in Section 19 (mostly written by A. Vladimirov) of \cite{KraPok}. More general cases of continuous convex moving sets $C(t)$ with non-empty interior were studied in \cite{Mar84, Mar86} (see also \cite[Chapter 2]{Mar93}). In all these references the fact that the set of constraints has non-empty interior allows the resulting solution to be of bounded variation.

All the above mentioned results deal with the case of convex constraints, however, the convexity assumption turns out to be too restrictive in some applications, for example, in problems coming from the modeling of crowd motion \cite{Venel}. The study of non-convex sweeping processes started with M. Valadier \cite{Val88} and, since then, has called the attention of many other authors, e.\,g., 
\cite{Ben00, ColGon99, Thi03}. An important concept which allows to get around the convexity of sets is the notion of uniform prox-regularity. These are closed sets having a neighborhood where the projection exists and is unique. Sets with such a property appear in the literature under different terminologies; being introduced under the name of `positively reached sets' by H. Federer \cite{Fed} in finite dimensional setting. A series of properties as well as the connection between sets and functions was deeply investigated in \cite{vial} (therein called `weak convex sets/functions').
The notion of prox-regularity was later extended to infinite dimensional spaces, \cite{clarke-95, prt}, and appears to lead to an appropriate class of non-convex sets for which one can prove existence and uniqueness results for sweeping processes, see for instance \cite{ANT, FV-12, FV-15, ChMM, EdT}. 
Notably, a recent paper \cite{NacThi} presents a fairly general result for $BV$ sweeping processes with prox-regular constraints. The case when the moving uniform prox-regular constraint has unbounded variation was instead dealt with in \cite{cmm} where it is assumed that $C(t)$ is continuous in time: In this paper another geometric condition is also required, namely $C(t)$ has uniform non-empty interior, which essentially means that cusps are not admitted on the boundary.

In the present paper we address the situation where the set of constraints is uniformly prox-regular and has uniform non-empty interior, but we also allow $C(t)$ to be discontinuous with possibly unbounded variation in time. We believe that the analysis of the problem becomes more transparent if in the motion of the set $C(t)$, we separate the effects of translation in the space $X$ from the effect of shape change, since in the mathematical description, translation and shape changes play completely different roles. To be more precise, we consider $C(t)$ of the form $C(t) = u(t) - Z(w(t))$ for given
$u : [0,T] \to X$ and $w(t) : [0,T] \to A$, where $A$ is a closed set of parameters in a Banach space $W$,
and we only assume that the translation $u$ and deformation $w$ are \emph{regulated right-continuous functions}, i.\,e., they admit the one-sided limits 
$u(t+)= u(t), u(t-), w(t+)=w(t), w(t-)$ at every point $t \in [0,T]$, with the convention $u(0-) = u(0)$, 
$w(0-) = w(0)$. Note that such functions are also called ``c\`adl\`ag" in the literature (= continue \`a droite, limite \`a gauche), see \cite{noe}. Concerning the shape of the moving constraint, we assume that $Z(w(t))$ is uniformly prox-regular and has uniformly non-empty interior. Indeed, since the projection onto a prox-regular constraint is defined only in a small neighborhood of the constraint, we have to keep the admissible jumps of the inputs $u$ and $w$ within suitable limits. 

The functional framework of regulated functions is convenient, since regulated functions are limits of step functions with respect to the topology of uniform convergence. We substantially make use of the Kurzweil integral calculus, which is compatible with the uniform convergence concept. Since Moreau's catching-up algorithm yields the exact solution in the Kurzweil integral setting for step functions $u$ and $w$, we obtain the general existence result in a standard way by passing to the limit.

The paper is structured as follows. In Section \ref{pro} we recall the notion and main properties of 
prox-regular set, and a rigorous statement of the main problem is specified in Section \ref{state}. In Section \ref{dis} we analyze a discretized version of our problem and derive a uniform bound for the output variation. Section \ref{reg} is devoted to the proof of convergence of the discrete scheme and of the continuous dependence property.
In Section \ref{con} we study the case when the inputs $u, w$ are continuous or absolutely continuous.
Finally in Appendix \ref{appe}, we collect some basic properties of the Kurzweil integral, which is a major tool in our analysis.


\section{Prox-regular sets}\label{pro}

\begin{definition}\label{d1}
Let $Z\subset X$ be a closed connected set and let $\dist(x,Z) := \inf\{|x-z|: z\in Z\}$ denote the distance of a point $x\in X$
from the set $Z$. Let $r>0$ be given. We say that $Z$ is {\em $r$-prox-regular} if the following condition hold.
\be{e1}
\forall y \in X:\ \dist(y,Z) = d \in(0, r) \ \ 
\exists\, x \in Z: \dist\left(x+\frac{r}{d}(y-x),Z\right) = \frac{r}{d}|y-x|= r.
\ee
\end{definition}

Note that this is in agreement with \cite[items (a) and (g) of Theorem 4.1]{prt}. We start with an easy lemma.

\begin{lemma}\label{l0}
Let $Z\subset X$ be an $r$-prox-regular set, and let $y\in X$ be given such that $\dist(y,Z) = d < r$.
Let $x$ satisfy the condition \eqref{e1}. For
$s \in [0,r]$ put $y(s) = x + (s/d)(y-x)$. Then $\dist(y(s), Z) =(s/d)|y-x| = s$ for every $s \in [0,r]$.
\end{lemma}

\bpf{Proof} For $s\in [0,r]$ we have $|y(s)- x|=s$. For every $z \in Z$ we have by the triangle inequality
$$
r \le |y(r)-z| \le \frac{r-s}{d}|y-x| + |y(s) - z| = r-s + |y(s) - z|, 
$$
hence, $|y(s) - z| \ge s$ for all $z \in Z$, which we wanted to prove.
\epf

For the reader's convenience, we explicitly state and prove an easy result going back to 
\cite[formula (1.2), (a) and (f) of Theorem 4.1]{prt}.

\begin{lemma}\label{l2}
A set $Z\subset X$ is $r$-prox-regular if and only if for every $y \in X$ such that $d =\dist(y,Z) < r$ there exists a unique $x\in Z$ such that $|y-x| = d$ and
\be{e2w}
\scal{y-x, x-z} + \frac{|y-x|}{2r} |x-z|^2 \ge 0 \quad \forall z \in Z.
\ee
\end{lemma}

\bpf{Proof} We first prove that for every $r$-prox-regular set $Z$ condition \eqref{e2w} holds. The case $d=0$ is trivial and it suffices to choose $x=y$. For $y \in X$ such that $\dist(y,Z) = d\in (0,r)$ we use \eqref{e1} and find $x \in Z$ such that $|y-x| = d$.
Put $\hat y = x + \frac {r}{d} (y-x)$. By \eqref{e1} we have $\dist(\hat y, Z) = r = |\hat y - x|$. Let now 
$z\in Z$ be arbitrary. We have
\be{e3w}
0 \le \frac12|\hat y - z|^2 - \frac12|\hat y - x|^2 = \scal{\hat y - x, x - z} + \frac12|x-z|^2 = 
\frac{r}{d} \left(\scal{y-x, x-z} + \frac{|y-x|}{2r} |x-z|^2\right)
\ee
and \eqref{e2w} follows.

To check that \eqref{e2w} holds for a unique $x\in Z$, assume that there exist $x_1, x_2 \in Z$ satisfying \eqref{e2w}. Then
\begin{align*}
\scal{y-x_1, x_1-x_2} + \frac{|y-x_1|}{2r} |x_1-x_2|^2 &\ge 0,\\
\scal{y-x_2, x_2-x_1} + \frac{|y-x_2|}{2r} |x_2-x_1|^2 &\ge 0.
\end{align*}
Summing up the above inequalities we obtain
$$
|x_1-x_2|^2 \le \frac{|y-x_1|+|y-x_2|}{2r}\, |x_1-x_2|^2 \le \frac{d}{r}\, |x_1-x_2|^2,
$$
hence $x_1 = x_2$, and the "only if" part of the proof is complete.

Assume now that for every $y \in X$ such that $d =\dist(y,Z) < r$ there exists a unique $x\in Z$ such that $|y-x| = d$ and \eqref{e2w} holds.
Let $y \in X$ be arbitrarily chosen such that $\dist(y,Z) = d \in (0,r)$. For $\hat y$ as in \eqref{e3w} we can now read \eqref{e3w} in the reverse
order and check that $|\hat{y}-z|\ge|\hat{y}-x|=\frac{r}{d}|y-x|$ for every $z\in Z$. Therefore 
$\dist(x+(r/d)(y-x),Z) = (r/d)|y-x|= r$, which we wanted to prove.
\epf

Clearly, every convex closed set $Z\subset X$ is $r$-prox-regular for all $r>0$. 
The vector $y-x$ in Lemma \ref{l2} is called {\em outward prox-normal vector\/}.
Indeed, an $r$-prox-regular set $Z$ admits a neighborhood $U_r(Z) := \{y \in X;\ \dist(y,Z) < r\}$ such that 
the mapping $P:U_r(Z) \to Z$ which with $y \in U_r(Z)$
associates $x\in Z$ from Lemma \ref{l2} is well defined and is called the {\em proximal projection\/} onto $Z$. Moreover, the set
\be{e2}
N_Z(x) = \{\xi \in X; \scal{\xi, x-z} + \frac{|\xi|}{2r} |x-z|^2 \ge 0 \quad \forall z \in Z\} 
\ee
is called the {\em proximal normal cone\/} to $Z$ at the point $x$, see, e.\,g., \cite{clarke-98}.

We always have $0 \in N_Z(x)$ for every $Z$ and every $x$. One might expect that the proximal normal cone to an $r$-prox-regular set contains nonzero elements for each $x \in \partial Z$.
Examples show that if $X$ is infinite-dimensional, a nonzero outward normal vector may fail to exist even in the convex case, and the existence is guaranteed if for example $Z$ is convex and
has non-empty interior (see, e.\,g., \cite[Proposition 2.11]{K-book}).
In the sequel, we therefore restrict the family of admissible sets $Z$ and we now
formulate a suitable non-empty interior condition for the nonconvex case. 
Here and in what follows for $x \in X$ and $\delta > 0$ we denote by $B_\delta(x)$ the closed ball
$\{y \in X;\ |y-x| \le \delta\}$.

\begin{definition}\label{d2}
Let $\mathcal{Z}$ be a family of $r$-prox-regular sets $Z\subset X$. We say that elements $Z\in \mathcal{Z}$ have {\em uniformly non-empty interior} if there exist $R\ge 3$ and $\rho \in (0, 2r/(1+R^2))$ such that for every $Z\in \mathcal{Z}$ we have
\be{e4}
\forall x \in Z \ \ \exists \bar x \in Z: \ |x-\bar x| < R\rho, \ B_{3\rho}(\bar x) \subset Z.
\ee
\end{definition}

Condition \eqref{e4} can be equivalently written as
\be{e4b}
\exists \rho \in \left(0,\frac{r}5\right)\ \ \forall x \in Z \ \ \exists \bar x \in Z: 
\ |x-\bar x|^2 +\rho^2 < 2r\rho, \ B_{3\rho}(\bar x) \subset Z.
\ee
Indeed, if \eqref{e4} holds, then $2r \rho > (1+R^2)\rho^2 > |x-\bar x|^2 + \rho^2$, which implies \eqref{e4b}. Conversely, if \eqref{e4b} holds, then there exists $R\ge 3$ such that
$$
\frac{|x-\bar x|^2}{\rho^2} < R^2 < \frac{2r}{\rho} - 1
$$
and \eqref{e4} follows. 

\begin{lemma}\label{l2a}
Let an $r$-prox-regular set $Z\subset X$ satisfy condition \eqref{e4} for some admissible values of 
$R>0, \rho > 0$. Then for every
$x \in \partial Z$ there exists a unit vector $\xi \in X$ such that 
\be{e3a}
\scal{\xi, x-z} + \frac{1}{2r} |x-z|^2 \ge 0 \quad \forall z \in Z.
\ee
\end{lemma}

\bpf{Proof}
Let $x \in \partial Z$ be given. We find a sequence $\{y_n; n\in \nat\}\subset X\setminus Z$ such that 
$\sup_{n\in\nat}\dist(y_n,Z)<r$ and $y_n \to x$ as $n \to \infty$. By \eqref{e1}
there exist $x_n\in X$ such that $\ve_n := |y_n - x_n| = \dist(y_n, Z) \le |y_n - x|$. From Lemma \ref{l2} it follows that the inequality
\be{e3b}
\scal{y_n - x_n, x_n-z} + \frac{|y_n - x_n|}{2r} |x_n-z|^2 \ge 0
\ee 
holds for every $n \in \nat$ and $z \in Z$. For all $n$ put
$$
\xi_n = \frac{1}{\ve_n}(y_n - x_n).
$$
Then
\be{e3d}
\scal{\xi_n, x_n-z} + \frac{1}{2r} |x_n-z|^2 \ge 0
\ee 
for every $n \in \nat$ and $z \in Z$.

We have $|x_n - x| \le |y_n - x| + |y_n - x_n|$, hence,
$$
\lim_{n\to \infty} |x_n - x| = 0.
$$
All $\xi_n$ are unit vectors. We can therefore find $\xi_\infty\in B_1(0)$ such that 
$\xi_n \rightharpoonup \xi_\infty$ weakly in~$X$.
In \eqref{e3d}, we choose
$$
z = \bar x + \rho \xi_n
$$
with the notation from \eqref{e4}, and obtain
$$
\scal{\xi_n, x_n-\bar x} - \rho + \frac{1}{2r} \left(|x_n-\bar x|^2 + \rho^2 - 2\rho\scal{\xi_n, x_n-\bar x} \right) \ge 0.
$$
Multiplying the above inequality by $2r$ we have by virtue of \eqref{e4} that
\be{e3c}
2(r-\rho) \scal{\xi_n, x_n-\bar x} \ge 2r\rho - \rho^2 - R^2\rho^2 
= \rho (2r - \rho(1+R^2)) =: \gamma >0.
\ee
Passing to the limit in \eqref{e3c} (note that $x_n$ converge strongly and $\xi_n$ weakly) we see that
\be{e3e}
2(r-\rho) \scal{\xi_\infty, x-\bar x} \ge \gamma >0. 
\ee
This implies, in particular, that $|\xi_\infty| = \sigma \in (0,1]$, and passing to the limit in \eqref{e3d} we conclude that
\be{e3f}
\scal{\xi_\infty, x-z} + \frac{1}{2r} |x-z|^2 \ge 0 \quad \forall z \in Z.
\ee 
This is precisely \eqref{e3a} for $\xi = \xi_\infty$ if $\sigma = 1$. For $\sigma<1$ we put $y = x + r\xi_\infty$. Then 
$$
\scal{y - x, x - z} + \frac12|x-z|^2 \ge 0 \quad \forall z \in Z,
$$
which is equivalent to the fact that $|y-z| \ge |y-x| = r\sigma$ for all $z \in Z$ (argue as in the computation in \eqref{e3w}). Thus, by Lemma \ref{l2} for $d= r\sigma$ we have
\be{e3g}
\scal{y-x, x-z} + \frac{|y-x|}{2r} |x-z|^2 \ge 0 \quad \forall z \in Z,
\ee
hence \eqref{e3a} holds for $\xi = \xi_\infty/\sigma$.
\epf

\begin{figure}[htb]
\begin{center}
\includegraphics[width=10cm]{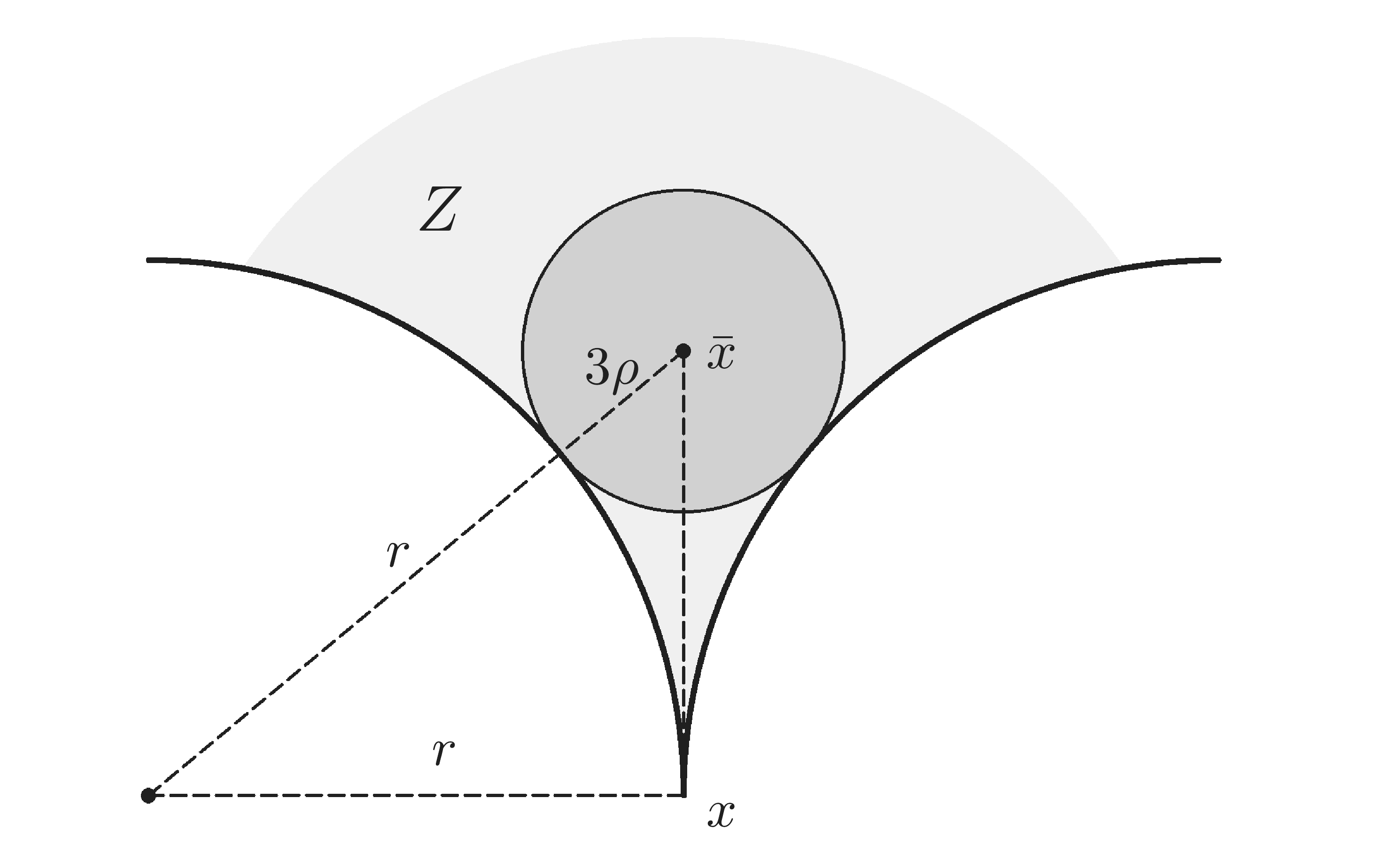}
\end{center}
\vspace{-9mm}

\caption{Violation of the uniform non-empty interior condition.}
\label{f1}
\end{figure}

In the situation of Figure \ref{f1}, where the set $Z$ admits a sharp cusp at $x$, for every $\rho > 0$ and every $\bar x$ positioned as in the picture we have  
that
$$
|x-\bar x|^2 +\rho^2 = 6r\rho + 10 \rho^2 > 2r\rho,
$$
and this is enough to show that condition \eqref{e4b} is violated.
In \cite[Theorem 4.2]{cmm}, the study of nonconvex sweeping processes relies on an interior cone condition which we write in the form
\be{p2}
\begin{aligned}
&\exists s>0\ \ \exists d>0 \ \ \forall x \in Z \ \ \exists x^* \in Z, |x-x^*| \le d : \\ 
&\ \forall \alpha \in (0,1], \ \ x + \alpha (x^* - x) + \alpha B_s(0) \subset Z. 
\end{aligned}
\ee
 Let us make the following interesting observation.


\begin{lemma}\label{lint}
Conditions \eqref{e4b} and \eqref{p2} are equivalent.
\end{lemma}

\bpf{Proof}
Assume first that \eqref{p2} holds. It suffices to put $\rho = \alpha s/3$ and $\bar x = x + \alpha (x^* - x)$. Then $|x-\bar x|^2 +\rho^2 - 2r\rho = \alpha^2(|x^* - x|^2 + s^2/9) - (2/3)\alpha r s
\le \alpha^2(d^2 + s^2/9) - (2/3)\alpha r s < 0$ for $\alpha$ sufficiently small.

Conversely, assume that \eqref{e4b} holds and suppose by contradiction that \eqref{p2} is not satisfied, so that in particular there exist $x \in Z$, 
$\alpha_0 \in (0,1]$ 
and $z \in B_1(0)$ with $z \neq 0$, such that 
$x_{\alpha_0} := x + \alpha_0 (\bar x - x + \rho z) \not\in Z$. Since 
$x_1 := x + (\bar x - x + \rho z) = \bar x + \rho z \in Z$ by virtue of \eqref{e4b},
we have that $\alpha_0 < 1$, and the segment connecting $x_{\alpha_0}$ and $x_1$
necessarily intersects the boundary $\partial Z$ of $Z$,
therefore there exists $\alpha \in (0,1)$ such that
$$
 x_\alpha := x + \alpha (\bar x - x) + \alpha \rho z \in \partial Z.
$$
Hence by Lemma \ref{l2a} and Lemma \ref{l2}, there exists $\xi \in X$, $|\xi| = 1$ such that 
$\dist(x_\alpha + r\xi, Z) = r$. By hypothesis, both $x$ and $\bar x +3\rho \xi$ belong to $Z$, hence
\begin{align*}
|x_\alpha + r\xi - x| \ge r,\\
|x_\alpha + (r-3\rho)\xi - \bar x| \ge r.
\end{align*}
In other words,
\begin{align*}
|\alpha(\bar x - x) +\alpha\rho z + r\xi| \ge r,\\
|(1-\alpha)(\bar x - x) - (r-3\rho)\xi - \alpha\rho z| \ge r,
\end{align*}
hence, using the triangle inequality and the fact that $|z| \le 1$,
\begin{align*}
|\alpha(\bar x - x) + r\xi| \ge r-\alpha\rho,\\
|(1-\alpha)(\bar x - x) - (r-3\rho)\xi|  \ge r -\alpha\rho,
\end{align*}
and, squaring both inequalities,
\begin{align*}
\alpha^2|\bar x - x|^2 + 2\alpha r\scal{\xi,\bar x - x} \ge -2\alpha r\rho +\alpha^2\rho^2,\\ 
(1-\alpha)^2|\bar x - x|^2 - 2(1-\alpha)(r-3\rho)\scal{\xi,\bar x - x}  \ge 
     6r\rho - 9\rho^2 -2\alpha r\rho + \alpha^2\rho^2.
\end{align*}
Taking into account that $r-3\rho > 0$ by \eqref{e4b}, we now eliminate the term 
$\scal{\xi,\bar x - x}$ from the above inequalities and obtain
\be{co1}
\begin{aligned}
(1-\alpha)\left(1-\alpha +\alpha \frac{r-3\rho}{r}\right)|\bar x - x|^2 &\ge r\rho \left(6 - 2\alpha - 2(1-\alpha)\frac{r-3\rho}{r}\right)\\
&\quad + \rho^2 \left(-9 +\alpha^2+\alpha(1-\alpha)\frac{r-3\rho}{r}\right).
\end{aligned}
\ee
Therefore, since $(r-3\rho)/r < 1$, from \eqref{e4b} and \eqref{co1} we infer that
$$
2r\rho - \rho^2 > |\bar x - x|^2 \ge 4r\rho - 9 \rho^2.
$$
This implies that $r<4\rho$, which contradicts \eqref{e4b}. The proof is complete.
\epf

In what follows, we assume that
\begin{equation}\label{W}
  \text{$W$ is a real Banach space endowed with norm $|\cdot|_W$}
\end{equation}
and that
\begin{equation}\label{A}
  \text{$A$ is a closed subset of $W$},
\end{equation}
and consider $r$-prox-regular sets $Z(w)$ depending on an additional parameter $w$ belonging
to $A$. We assume that  the dependence of $Z$ on $w$ is continuous with respect to the Hausdorff distance
$$
d_H(Z,\hat Z) := \max\{\sup_{z\in Z}\dist(z, \hat Z), \sup_{\hat z\in \hat Z}\dist(\hat z, Z)\},
$$
more specifically, we assume that
\be{ha1}
\forall \ve > 0 \ \exists \delta >0: \ |w-\hat w|_W < \delta \ \Longrightarrow \ d_H(Z(w), Z(\hat w)) < \ve.
\ee


\section{Statement of the problem}\label{state}

In this section we provide a rigorous formulation of the non-convex sweeping process in the framework of regulated functions. We first recall some basic facts about regulated functions which can be found, e.\,g., in \cite{H}.

\begin{definition}\label{s1}
Assume that \eqref{W} holds. A function $f : [0,T] \to W$ is called \emph{regulated (on $[0,T]$)} if it admits the one sided limits $f(t-)$, $f(t+)$ for every $t \in [0,T]$, with the convention that $f(0-) := f(0)$ and 
$f(T+) := f(T)$. The space of $W$-valued regulated functions on $[0,T]$ is denoted by $G(0,T;W)$ and we also set $G_R(0,T;W) := \{f \in G(0,T;W);\ f(t+) = f(t)\ \forall t \in [0,T]\}$. 
\end{definition}

Indeed, every regulated function is bounded and the set of its discontinuity points is at most countable, cf. \cite[Corollary I.3.2]{H}. Moreover, the space $G(0,T;W)$ endowed with the supremum norm
\be{nw}
\|f\| := \sup_{t \in [0,T]} |f(t)|_W \ \for f \in G(0,T;W)
\ee
is a Banach space and $G_R(0,T;W)$ is its closed subspace.

An important dense subset of $G(0,T;W)$ consists of the so-called \emph{step functions}, that is, functions $f : [0,T] \to W$ such that there is a division $0 = t_0 < t_1 < \cdots < t_m = T$ for which $f$ is constant on any interval $(t_{j-1},t_j)$, $j=1,\ldots,m$.

Another proper subset of $G(0,T;W)$ important in our investigation is provided by the space $BV(0,T;W)$ of functions of bounded variation which we now briefly recall. For $f : [0,T] \to W$ and for $a,b \in [0,T]$, $a < b$, we set
\begin{equation}\label{s3}
  \Var_{[a,b]} f := 
  \sup\left\{\sum_{j=1}^m|f(t_j) - f(t_{j-1})|_W;\ a = t_0 < \cdots < t_m = b,\ m \in \nat\right\},
\end{equation}
and we define $BV(0,T;W) := \{f : [0,T] \to W;\ \Var_{[0,T]} f < \infty$\} and
$BV_R(0,T;W) := BV(0,T;W) \cap G_R(0,T;W)$. It is also convenient to set 
\begin{equation}\label{s3-2}
  V(f)(t) := \Var_{[0,t]} f
\end{equation} 
for $t \in [0,T]$ and $f \in BV(0,T;W)$. The function $V(f)$ is bounded and nondecreasing for every 
$f \in BV(0,T;W)$, hence $V(f) \in BV(0,T;\real)$.

The following result will be used throughout
the paper (see, e.\,g., \cite[Theorem I.3.1]{H} or \cite{Fra}).

\begin{proposition}\label{s2}
For each $f \in G(0,T;W)$ there exists a sequence of step functions $f_k : [0,T] \to W$, $k \in \nat$ which is uniformly convergent to $f$ on $[0,T]$ and such that $f_k([0,T]) \subset f([0,T])$. In particular, if 
$f \in G_R(0,T;W)$ one can assume that every $f_k$ is right continuous, and if $f \in BV(0,T;W)$ then 
$\Var_{[0,t]} f_k \le \Var_{[0,t]} f$ for every $k\in \nat$ and $t \in [0,T]$.
\end{proposition} 

Now we are ready to state the precise formulation of our main problem.

\begin{problem}\label{PB}
Assume that \eqref{i1} and \eqref{W}--\eqref{A} hold and that $r > 0$. Let $\{Z(w);\ w \in A\}$ be a family of $r$-prox-regular subsets of $X$
satisfying \eqref{ha1}. For given functions $u \in G_R(0,T;X)$, $w \in G_R(0,T;W)$ such that $w(t) \in A$ for every $t \in [0,T]$ and $x_0 \in Z(w(0))$, we look for functions $\xi \in BV_R(0,T;X)$ and 
$x \in G_R(0,T;W)$ such that $\xi(t) + x(t) = u(t)$ for all $t \in [0,T]$ and 
\begin{align}
  & \int_0^T \scal{x(t) - z(t), \dd\xi(t)} + \frac1{2r}\int_0^T |x(t) - z(t)|^2 \dd V(\xi)(t) \ge 0 \notag \\
  & \qquad\qquad\qquad\quad \forall z \in G(0,T;X), \  \text{$z(t) \in Z(w(t))$ for all $t \in [0,T]$,} \label{s4}\\
  & x(0) = x_0. \label{s5}
\end{align}
\end{problem}

The two integrals in \eqref{s4} are meant in the sense of the Kurzweil integral introduced in \cite{kur57}: In the first integral we are integrating $X$-valued functions, while on the second integral the particular case of real-valued functions is considered.
In Appendix \ref{appe} we briefly recall the main elements of the Kurzweil integral calculus. A comparison with \eqref{e2} shows that \eqref{s4} can be interpreted as an integral formulation of the inclusion
\eqref{i1a} with $C(t) = u(t) - Z(w(t))$.
The solution of Problem \ref{PB} will be provided by Theorem \ref{t1} below under an additional assumption \eqref{e5g}.


\section{Discrete sweeping processes}\label{dis}

In this section we study a family of discrete sweeping processes which are obtained by an implicit Euler discretization scheme for inequality \eqref{s4}.

\begin{lemma}\label{l1}
We assume that $\{Z(w); w \in A\subset W\}$ is a family of $r$-prox-regular subsets of $X$ satisfying condition \eqref{ha1}, and that $\{u_j; j\in\nat\cup\{0\}\} \subset X$, $\{w_j; j\in\nat\cup\{0\}\}\subset A$ are given sequences such that
\be{e5}
d_H(Z(w_j), Z(w_{j-1})) + |u_j - u_{j-1}| \le \frac{r}{M} \quad \forall j\in\nat
\ee
for some $M\ge 2$. We further assume that $x_0 \in Z(w_0)$ is a given element and put $\xi_0 := u_0 - x_0$. Then there exists a unique pair of sequences $\{\xi_j; j\in\nat\}\subset X$,  
$\{x_j; j\in\nat\} \subset X$ such that 
\be{xxi}
x_j \in Z(w_j), \ \xi_j = u_j - x_j, \ \scal{\xi_j - \xi_{j-1}, x_j - z} + \frac{|\xi_j - \xi_{j-1}|}{2r}|x_j - z|^2 \ge 0 \ \ \forall z \in Z(w_j).
\ee
Moreover the following estimates hold  for every $j \in \nat$:
\be{e7}
\begin{aligned}
|x_j - x_{j-1}| &\le \frac{2M}{2M-1} |u_j - u_{j-1}|+ \frac{4M-1}{2M-1} d_H(Z(w_j), Z(w_{j-1})), \\[3mm]
|\xi_j - \xi_{j-1}| &\le \frac{4M-1}{2M-1}\big(|u_j - u_{j-1}|+  d_H(Z(w_j), Z(w_{j-1}))\big).
\end{aligned}
\ee
\end{lemma}

The exact value of the constant $M$ is not important for the moment. It will be specified below in Theorem \ref{t1}.

\bpf{Proof of Lemma \ref{l1}}
The sequences $x_j$ and $\xi_j$ can be uniquely constructed in the following recursive way. Assume that we have already constructed $x_i \in Z(w_i)$ and $\xi_i = u_i - x_i$ for $i=0, \dots, j-1$, and put 
$y = x_{j-1} + u_j - u_{j-1}$. From \eqref{e5} it follows that $d = \dist(y,Z(w_{j})) \le r/M$, thus by 
Lemma \ref{l2} there exists a unique $x \in Z(w_{j})$ such that $|y-x| = d$ and
$$
\scal{y-x, x-z} + \frac{|y-x|}{2r} |x-z|^2 \ge 0 \quad \forall z \in Z(w_j).
$$
Putting $x_j := x$ and $\xi_j := u_j - x_j$ we obtain $y-x = \xi_j - \xi_{j-1}$ and \eqref{xxi} follows. Now by construction we have
\be{e6}
|\xi_j - \xi_{j-1}| \le \frac{r}{M}
\ee
for every $j \in \nat$, and Lemma \ref{l2} enables us to find $\hat x_{j-1} \in Z(w_j)$ such that 
$|x_{j-1} - \hat x_{j-1}| = \dist(x_{j-1}, Z(w_j))$. Therefore putting $z=\hat x_{j-1}$ in \eqref{xxi} yields
\be{e6+}
\scal{\xi_j - \xi_{j-1}, x_j - \hat x_{j-1}} + \frac{|\xi_j - \xi_{j-1}|}{2r}|x_j - \hat x_{j-1}|^2 \ge 0\ \for j\in \nat,
\ee
which implies
$$
\begin{aligned}
  |x_j - \hat x_{j-1}|^2 
    & \le \scal{u_j - u_{j-1}, x_j - \hat x_{j-1}} + \frac{|\xi_j - \xi_{j-1}|}{2r}|x_j - \hat x_{j-1}|^2
            + \scal{x_{j-1} - \hat x_{j-1}, x_j - \hat x_{j-1}} \\ 
    & \le |u_j - u_{j-1}|\,|x_j - \hat x_{j-1}| + \frac1{2M}|x_j - \hat x_{j-1}|^2
            + |x_{j-1} - \hat x_{j-1}|\,|x_j - \hat x_{j-1}|.
\end{aligned}
$$
Hence,
$$
\left(1 - \frac1{2M}\right)|x_j - \hat x_{j-1}|\le d_H(Z(w_j), Z(w_{j-1})) + |u_j - u_{j-1}|,
$$
which implies that
\begin{align*}
|x_j - x_{j-1}| &\le |x_j - \hat x_{j-1}|+|x_{j-1} - \hat x_{j-1}|\\
&\le
\frac{2M}{2M-1} \big(d_H(Z(w_j), Z(w_{j-1})) + |u_j - u_{j-1}|\big) + d_H(Z(w_j), Z(w_{j-1}))
\end{align*}
and we easily obtain the upper bounds \eqref{e7}.
\epf


\subsection{Estimates of the total variation}

Let the assumptions of Lemma \ref{l1} 
hold and let $\{\xi_j\}$ be the sequence defined by \eqref{xxi}. It is easy to estimate the output variation 
$\sum_{j=1}^n |\xi_j - \xi_{j-1}|$ for any $n\in \nat$ using \eqref{e7} if we control the input variation 
$\sum_{j=1}^n \big(|u_j - u_{j-1}|+ d_H(Z(w_j), Z(w_{j-1}))\big)$. In this subsection, we show that if $Z(w)$ have uniformly non-empty interior, the variation of $\{\xi_j\}$ can be estimated even if the variations of 
$\{u_j\}$ and $\{w_j\}$ are unbounded. The argument is based on the following statement.

\begin{proposition}\label{p1}
Let $\{Z(w); w \in A\}$, $\{u_j\}, \{w_j\}$, and $M$ satisfy the assumptions of Lemma \ref{l1}.  Suppose further that all elements of the system $\{Z(w); w \in A\}$ have uniformly non-empty interior with $\rho, R$ as in Definition \ref{d2}, and assume that there exist $j_0 < j_1$ such that
\be{e8}
|u_j - u_{j_0}| + d_H(Z(w_j), Z(w_{j_0})) \le \rho \quad \for j= j_0+1, \dots, j_1.
\ee
Then, if $\{x_j\}, \{\xi_j\}$ satisfy \eqref{xxi} and putting $\mu := \rho(2r - \rho)$, we have that
\be{e9}
\sum_{j = j_0+1}^{j_1} |\xi_j - \xi_{j-1}| \le (r-\rho)\log\left(\frac{\mu}{\mu - R^2\rho^2}\right),
\ee
where the right-hand side of \eqref{e9} makes sense since by Definition \ref{d2} we have 
$(1+R^2)\rho < 2r$ so that 
\be{e9b}
R^2\rho^2 < \mu.
\ee
\end{proposition}

\bpf{Proof}
In \eqref{xxi} it suffices to consider those values of $j$ for which $|\xi_j - \xi_{j-1}| > 0$. By Definition \ref{d2} we choose $\bar x \in Z(w_{j_0})$ such that
\be{e9a}
|u_{j_0} - \xi_{j_0} - \bar x| \le R\rho, \quad B_{3\rho}(\bar x) \subset Z(w_{j_0}),
\ee
and put in \eqref{xxi} for $j= j_0+1, \dots, j_1$
\be{test}
z = u_j - u_{j_0} + \rho \frac{\xi_j - \xi_{j-1}}{|\xi_j - \xi_{j-1}|} + \bar x.
\ee
This is an admissible choice provided we show that
\be{2rho}
B_{2\rho}(\bar x) \subset Z(w_j)\  \for j= j_0+1, \dots, j_1.
\ee
Indeed, assume that there exists $j \in \{j_0+1, \dots, j_1\}$ and $y\in X$ such that $|y-\bar x| \le 2\rho$ and $y \notin Z(w_j)$.
We have $y\in Z(w_{j_0})$, hence $\dist(y,Z(w_{j})) =:\delta_j \in (0,\rho]$. Let $\bar x_j \in Z(w_{j})$ be such that $|y - \bar x_j| = \delta_j$.
For $s \in [0,r]$ put
$$
y_j(s) = \bar x_j + \frac{s}{\delta_j}(y - \bar x_j).
$$
By Lemma \ref{l0} we have $\dist(y_j(s), Z(w_j)) = |y_j(s) - \bar x_j| = s$. On the other hand,
for $\rho< s \le \rho+\delta_j$ we have
$$
|y_j(s) - \bar x| \le |y-\bar x| + \left(\frac{s}{\delta_j} - 1\right)|y - \bar x_j|  \le 3\rho,
$$ 
hence $y_j(s) \in Z(w_{j_0})$, $\dist(y_j(s), Z(w_j)) = |y_j(s) - \bar x_j| = s >\rho$, which is a contradiction, and \eqref{2rho} is proved.
It follows that $z$ given by \eqref{test} can be chosen as test element in \eqref{xxi} for all $j= j_0+1, \dots, j_1$.
We obtain
\begin{align*}
&\scal{\xi_j - \xi_{j-1}, u_{j_0} - \xi_j - \bar x} - \rho|\xi_j - \xi_{j-1}|\\
&\quad + \frac{|\xi_j - \xi_{j-1}|}{2r}\left(|u_{j_0} - \xi_j - \bar x|^2 + \rho^2 - \frac{2\rho}{|\xi_j - \xi_{j-1}|}\scal{\xi_j - \xi_{j-1}, u_{j_0} - \xi_j - \bar x}\right) \ge 0,
\end{align*}
that is
\be{e10}
\left(1 - \frac{\rho}{r}\right)\scal{\xi_j - \xi_{j-1}, u_{j_0} - \xi_j - \bar x} + \frac{|\xi_j - \xi_{j-1}|}{2r}|u_{j_0} - \xi_j - \bar x|^2
\ge \rho\left(1 - \frac{\rho}{2r}\right) |\xi_j - \xi_{j-1}|.
\ee
Put
\be{e11}
U_j = |u_{j_0} - \xi_j - \bar x|^2,\quad s_j = \frac{|\xi_j - \xi_{j-1}|}{r-\rho}.
\ee
Using the elementary inequality $\scal{a, a-b} \ge \frac12 (|a|^2 - |b|^2)$ for $a,b \in X$ we easily conclude that
\be{e12}
\scal{\xi_j - \xi_{j-1}, u_{j_0} - \xi_j - \bar x} \le -\frac12 (U_j - U_{j-1}).
\ee
Combining \eqref{e10}, \eqref{e11}, and \eqref{e12} we get
\be{e13}
U_{j-1} - U_j + s_j U_j \ge \mu s_j
\ee
with $\mu = \rho(2r - \rho)$ as in the hypotheses.
Inequality \eqref{e13} holds trivially if  $\xi_j = \xi_{j-1}$, so that it is fulfilled for all
$j= j_0+1, \dots, j_1$.

We first rewrite \eqref{e13} as
\be{e14}
(1-s_j) U_j \le U_{j-1} - \mu s_j.
\ee
We have $M\ge 2$ and $\rho < r/5$. Hence, from \eqref{e6} and \eqref{e4b} it follows that $s_j \le r/M(r-\rho) < 5/8 < 1$. Furthermore, \eqref{e9b}--\eqref{e9a} yield that
\be{e14a}
U_{j_0} \le R^2 \rho^2 < \mu.
\ee
By induction it follows from \eqref{e14} that $U_j < \mu$ for all $j= j_0+1, \dots, j_1$. Using \eqref{e13} we conclude that $U_j \le U_{j-1}$ for all $j= j_0+1, \dots, j_1$, and we rewrite \eqref{e13} again as
\be{e15}
s_j \le \frac{U_{j-1} - U_j}{\mu - U_j}.
\ee
For all $v \in [U_j, U_{j-1}]$ we have $\mu - U_j \ge \mu - v$, hence
$$
s_j \le \frac{U_{j-1} - U_j}{\mu - U_j} \le \int_{U_j}^{U_{j-1}}\frac{\dd v}{\mu - v} = \log(\mu - U_j) - \log(\mu- U_{j-1}).
$$
Summing up the above inequalities over $j$ we obtain
\be{e16}
\sum_{j = j_0+1}^{j_1} |\xi_j - \xi_{j-1}| \le (r-\rho)\log\left(\frac{\mu- U_{j_1}}{\mu - U_{j_0}}\right)
\ee
and the assertion follows from \eqref{e14a}.
\epf

Note that \eqref{xxi} is the so-called {\em catching up algorithm\/}, see \cite{moreau}.  
It can also be interpreted in terms of the Kurzweil integral of piecewise constant functions, and we explain this approach in the next Subsection.

The uniform nonempty interior condition \eqref{e4b} is necessary for the validity of the upper bound \eqref{e16} if the input variation is unbounded. In the situation of Figure \ref{f1}, we can consider a fixed set $Z(w) = Z$, the initial condition $x_0 = x$ at the tip of the cusp, and an input $\{u_j\}$ which oscillates perpendicularly to the vector $\bar x - x$. Then we easily check that the recipe \eqref{xxi} yields $\xi_j = u_j$. The variation of $\{\xi_j\}$ is
therefore the same as the variation of $\{u_j\}$ which can be arbitrarily large.


\subsection{Kurzweil integral sweeping processes}\label{kur}

Consider right continuous step functions
$u: [0,T] \to X$, $w: [0,T] \to A\subset W$ of the form
\be{k1}
\begin{aligned}
u(t) &= \sum_{j=1}^{m} u_{j-1} \chi_{[t_{j-1}, t_j)}(t) + u_m \chi_{\{t_m\}}(t),\\ 
w(t) &= \sum_{j=1}^{m} w_{j-1} \chi_{[t_{j-1}, t_j)}(t) + w_m \chi_{\{t_m\}}(t), 
\end{aligned}
\ee
corresponding to a division $0 = t_0 < t_1 < \dots < t_m = T$ of the interval $[0,T]$. 
In \eqref{k1}, $\chi_A$ denotes the characteristic function of $A \subset [0,T]$, that is, $\chi_A(t) = 1$ for $t \in A$, $\chi_A(t) = 0$ for $t \in [0,T]\setminus A$. 

Let $u_j, w_j$ satisfy the hypotheses of Lemma \ref{l1}, let $x_j, \xi_j$ be associated with $u_j$, $w_j$ as in \eqref{xxi}, and put
\begin{align}\label{k2}
x(t) &= \sum_{j=1}^{m} x_{j-1} \chi_{[t_{j-1}, t_j)}(t) + x_m \chi_{\{t_m\}}(t),\\ \label{k3} 
\xi(t) &= \sum_{j=1}^{m} \xi_{j-1} \chi_{[t_{j-1}, t_j)}(t) + \xi_m \chi_{\{t_m\}}(t). 
\end{align}
Then $x(t) \in Z(w(t))$ for all $t\in [0,T]$, and \eqref{xxi} can be written in Kurzweil integral form \eqref{s4}, indeed, on one hand, by Theorems \ref{add-int} and \ref{cases}
we have
\begin{align}
  \int_0^T \scal{x(t) - z(t), \dd\xi(t)}  
  & = \sum_{j=1}^{m} \int_{t_{j-1}}^{t_j} \scal{x(t) - z(t), \dd(\xi_{j-1}\chi_{[t_{j-1}, t_j)}+\xi_{j} \chi_{\{t_{j}\}})(t)}  	\notag \\
     & = \sum_{j=1}^{m} \scal{x_j - z(t_j), \xi_{j} - \xi_{j-1}}. \label{k4lefthand}
\end{align}
On the other hand, observing that
$$
V(\xi)(t) = \sum_{j=1}^{m-1}(\sum_{k=1}^{j} |\xi_k - \xi_{k-1}|) \chi_{[t_{j}, t_{j+1})}(t) + 
\sum_{k=1}^{m} |\xi_k - \xi_{k-1}|\chi_{\{t_m\}}(t)
$$
and using Theorem \ref{cases} we get
$$
\int_{t_{j-1}}^{t_{j}} |x(t) - z(t)|^2 \dd V(\xi)(t)=|x_{j} - z(t_{j})|^2 |\xi_{j} - \xi_{j-1}|\quad\for j=1,\dots,m,
$$ 
and consequently
\begin{align}
  \int_0^T |x(t) - z(t)|^2 \dd V(\xi)(t) 
  & =  \sum_{j=1}^{m} \int_{t_{j-1}}^{t_j} |x(t) - z(t)|^2 \dd V(\xi)(t) \notag \\
  & = \sum_{j=1}^{m} |x_j - z(t_j)|^2 |\xi_j - \xi_{j-1}|. \label{k4righthand}
\end{align}
Thus it follows from \eqref{k4lefthand}--\eqref{k4righthand} that \eqref{xxi} and \eqref{s4} are equivalent.

The following property of Kurzweil integral variational inequalities will be useful in the sequel.

\begin{lemma}\label{kl1}
Let $x\in G_R(0,T;X)$ and $\xi \in BV_R(0,T;X)$  
satisfy \eqref{s4}. Then for every $0\le \sigma < \tau \le T$ we have
\be{k5}
\int_\sigma^\tau \scal{x(t) - \tilde z(t), \dd\xi(t)} + \frac1{2r}\int_\sigma^\tau |x(t) - \tilde z(t)|^2 \dd V(\xi)(t) \ge 0
\ee
for every regulated test function $\tilde z: [\sigma,\tau] \to X$, $\tilde z(t) \in Z(w(t))$ for all $t\in [\sigma,\tau]$.
\end{lemma}

\bpf{Proof} 
The argument is standard and follows the lines of the proof of \cite[Lemma 2.2]{kuhys}. It consists in choosing $z(t)$ in \eqref{s4} in the form
$$
z(t) = \left\{
\begin{array}{ll}
x(t) \for t\in [0,\sigma]\cup (\tau,T],\\
\tilde z(t) \for t \in (\sigma,\tau].
\end{array}
\right.
$$
\epf


\section{Arbitrary right continuous regulated inputs}
\label{reg}

One of the main tools in further analysis of the Kurzweil integral variational inequality \eqref{k5} is the Kurzweil integral counterpart of the Gronwall Lemma which goes back to \cite[Chapter~22]{kur}. For the reader's convenience, 
we prove Lemma \ref{gl2} in Appendix \ref{appe} as a simplified version of the general theory which is sufficient for our purposes. We start by proving a  general uniqueness result.

\begin{lemma}\label{unique}
Let $\{Z(w); w\in A\} \subset X$ be a family of $r$-prox-regular sets satisfying \eqref{ha1} and let 
$u \in G_R(0,T;X)$, $w\in G_R(0,T;W)$, $w(t) \in A$ for all $t \in [0,T]$, and $x_0 \in Z(w(0))$
be given. Then there exists at most one function $\xi \in BV_R(0,T;X)$ such that \eqref{s4}--\eqref{s5} hold with $x(t) = u(t) - \xi(t)$ for every $t \in [0,T]$.
\end{lemma}

\bpf{Proof}
Assume by contradiction that there exist two solutions $\xi_1, \xi_2 \in BV_R(0,T;X)$ of the variational inequality
\be{c11}
\int_0^T \scal{u(t) - \xi_i(t) - z(t), \dd\xi_i(t)} + \frac1{2r}\int_0^T |u(t) - \xi_i(t) - z(t)|^2 \dd V(\xi_i)(t) \ge 0
\ee
for every $z \in G(0,T;X)$, $z(t) \in Z(w(t))$ for all $t \in [0,T]$, $i=1,2$. By Lemma \ref{kl1}, we have for all $\tau \in [0,T]$ that
\be{c12}
\int_0^\tau \scal{u(t) - \xi_i(t) - \tilde z(t), \dd\xi_i(t)} + \frac1{2r}\int_0^\tau |u(t) - \xi_i(t) - \tilde z(t)|^2 \dd V(\xi_i)(t) \ge 0
\ee
for every admissible $\tilde z$, $i=1,2$. In the variational inequality \eqref{c11} we now choose $\tilde z(t) = u(t) - \xi_{2}(t)$ for $i=1$ and $\tilde z(t) = u(t) - \xi_{1}(t)$ \for $i=2$,
and sum the two inequalities up. This yields for $\bar \xi = \xi_1 - \xi_2$ that
\be{c13}
\int_0^\tau \scal{\bar\xi(t), \dd\bar\xi(t)} \le \frac1{2r}\int_0^\tau |\bar\xi(t)|^2 \dd (V(\xi_1)+V(\xi_2))(t).
\ee
We have indeed $\bar\xi(0) = 0$, and $\int_0^\tau \scal{\bar\xi(t), \dd\bar\xi(t)} \ge \frac12 |\bar\xi(\tau)|^2$ by Corollary \ref{3c17}. We are thus in the situation of Lemma \ref{gl2} which yields $\bar \xi \equiv 0$, hence, $\xi_1 = \xi_2$ and the proof is complete.
\epf

We now state the main existence theorem for \eqref{s4} for the case of right-continuous regulated inputs with moderate jumps under two different hypotheses, namely that either the sets $Z(w)$ have 
uniformly non-empty interior condition, or, alternatively, the inputs have bounded variation.

\begin{theorem}\label{t1}
Let $\{Z(w); w\in A\} \subset X$ be a family of $r$-prox-regular sets satisfying \eqref{ha1}
and let $u \in G_R(0,T;X)$, $w\in G_R(0,T;W)$, $w(t) \in A$ for all $t \in [0,T]$ be given such that
\be{e5g}
r^* := \sup_{t\in (0,T]} \left(d_H(Z(w(t)), Z(w(t-))) + |u(t) - u(t-)|\right) < \frac{r}{M},\quad M= \frac{9 + \sqrt{65}}{4}.
\ee
Assume furthermore that at least one of the following two conditions holds:
\begin{itemize}
\item[{\rm (i)}] 
The sets $\{Z(w); w\in A\}$ have uniformly non-empty interior with $\rho, R$ as in Definition \ref{d2}.
\item[{\rm (ii)}] $A=W$, $u \in BV_R(0,T;X)$, $w\in BV_R(0,T;W)$, and 
\be{ha2}
\exists L>0 \ \ \forall w_1, w_2 \in W: d_H(Z(w_1),Z(w_2)) \le  L |w_1 - w_2|_W.
\ee
\end{itemize}
Then for every initial condition $x_0 \in Z(w(0))$ there exists a unique $\xi \in BV_R(0,T;X)$
and a unique $x \in G_R(0,T;X)$, $x(t) \in Z(w(t))$, $\xi(t) + x(t) = u(t)$ for all $t\in [0,T]$, such that \eqref{s4} holds and $x(0) = x_0$. Moreover, the piecewise constant approximations $\xi\ok$ defined in \eqref{c5} below and associated with the catching up algorithm \eqref{xxi} converge uniformly to~$\xi$.
\end{theorem}

\bpf{Proof}
First of all let us consider some ${\nu} \le r^*$ and a division $0 = \hat t_0 < \hat t_1 < \dots < \hat t_N = T$ such that
\be{c1}
d_H(Z(w(t)), Z(w(\hat t_{i-1})))+|u(t) - u(\hat t_{i-1})| \le {\nu}
\ \for t \in [\hat t_{i-1}, \hat t_i), 
\ i=1, \dots, N.
\ee
The idea of \eqref{c1} is to isolate the points where the jumps are higher than ${\nu}$, and estimate the variation of $\xi$ separately inside the intervals $[\hat t_{i-1}, \hat t_i)$, and across the big jumps. The number $N$ depends on ${\nu}$, indeed, but for a regulated function, it is always finite.

Let $\{u\ok\}, \{w\ok\}$ be sequences 
of right continuous step functions which converge uniformly on $[0,T]$ to $u$, $w$, respectively, as $k \to \infty$, and such that
\be{e5gk}
\sup_{t\in (0,T]} \left(d_H(Z(w\ok(t)), Z(w\ok(t-))) + |u\ok(t) - u\ok(t-)|\right) \le r^*
\ee
for all $k \in \nat$. Thanks to \eqref{ha1} we find $k_0 \in \nat$ such that for $k \ge k_0$ we have
\be{c2}
\sup_{t\in [0,T]}d_H(Z(w(t)), Z(w\ok(t))) + \| u - u\ok\| < \frac{{\nu}}{2}.
\ee
We now choose arbitrary $k, l \ge k_0$. We can assume that the functions $u\ok, u\ol, w\ok, w\ol$ are step functions of the form
\be{c4}
\begin{aligned}
u\ok(t) &= \sum_{j=1}^{m} u\ok_{j-1} \chi_{[t_{j-1}, t_j)}(t) + u\ok_m \chi_{\{t_m\}}(t),\\  
u\ol(t) &= \sum_{j=1}^{m} u\ol_{j-1} \chi_{[t_{j-1}, t_j)}(t) + u\ol_m \chi_{\{t_m\}}(t), 
\end{aligned}
\ee
\be{c4w}
\begin{aligned}
w\ok(t) &= \sum_{j=1}^{m} w\ok_{j-1} \chi_{[t_{j-1}, t_j)}(t) + w\ok_m \chi_{\{t_m\}}(t),\\  
w\ol(t) &= \sum_{j=1}^{m} w\ol_{j-1} \chi_{[t_{j-1}, t_j)}(t) + w\ol_m \chi_{\{t_m\}}(t), 
\end{aligned}
\ee
where $0 = t_0 < t_1< \dots < t_m = T$ is a division containing all discontinuity points of all functions 
$u\ok, u\ol, w\ok, w\ol$ as well as the points $\hat t_0, \dots, \hat t_N$; moreover, we assume $u\ok_0=u\ol_0=u_0$. By \eqref{e5gk} we have
\begin{equation}\label{e5-k}
d_H(Z(w_j\ok), Z(w_{j-1}\ok)) + |u_j\ok - u_{j-1}\ok| \le r^* < r
\end{equation}
for $j=1,\dots\,m$.
Hence condition \eqref{e5} holds
with $w_j$, $u_j$ replaced with $w\ok_{j}$, $u\ok_{j}$ etc., and the corresponding solutions of \eqref{s4} for $k,l\ge k_0$ 
from Lemma \ref{l1} have the form
\be{c5}
\begin{aligned}
\xi\ok(t) &= \sum_{j=1}^{m} \xi\ok_{j-1} \chi_{[t_{j-1}, t_j)}(t) + \xi\ok_m \chi_{\{t_m\}}(t),\\ 
\xi\ol(t) &= \sum_{j=1}^{m} \xi\ol_{j-1} \chi_{[t_{j-1}, t_j)}(t) + \xi\ol_m \chi_{\{t_m\}}(t), 
\end{aligned}
\ee
where $\xi\ok_{j}$ (respectively $\xi\ol_{j}$) satisfy \eqref{xxi} with $u_j$, $w_j$, $x_j$, $\xi_j$ replaced with $u\ok_j$, $w\ok_j$, $x\ok_j$, $\xi\ok_j$  (respectively $u\ol_j$, $w\ol_j$, $x\ol_j$, $\xi\ol_j$), and where
$x_0^{(l)}=x\ok_0 = x_0$, $x\ok(t) = u\ok(t) - \xi\ok(t)$, $x\ol(t) = u\ol(t) - \xi\ol(t)$. 

We now make use of hypothesis \eqref{c1} and of the triangle inequality to obtain for all $k\ge k_0$ and $t \in [\hat t_{i-1}, \hat t_i)$, $i=1, \dots, N$ that
\be{c1k}
d_H(Z(w\ok(t)), Z(w\ok(\hat t_{i-1})))+|u\ok(t) - u\ok(\hat t_{i-1})| \le 2{\nu}.
\ee
We now distinguish between the two cases (i) and (ii) specified in Theorem \ref{t1}. In Case (i), we choose
\be{sig1}
{\nu} = \min\left\{r^*, \frac{\rho}{2}\right\}
\ee
with $\rho$ from Definition \ref{d2}. Then condition \eqref{e8} is satisfied, so that we can use Proposition \ref{p1} for estimating the output variation. Indeed, We have
\be{vark}
\Var_{[0,T]} \xi\ok = \sum_{i=1}^N \left(\Var_{[\hat t_{i-1},\hat t_i)} \xi\ok + |\xi\ok(\hat t_i) - \xi\ok (\hat t_i-)|\right) 
\ee
and similarly for $\xi\ol$. 
Thus, we get from Proposition \ref{p1}, \eqref{e5gk}, and \eqref{e7} the upper bound
\be{c6}
\max\left\{\Var_{[0,T]} \xi\ok, \Var_{[0,T]} \xi\ol\right\} \le N(r-\rho)\log\left(\frac{\mu}{\mu - R^2\rho^2}\right) + 3Nr^* =: \bar C,
\ee
which is independent of $k$ and $l$ and depends on $u$ only 
through the number $N$.
In Case (ii), we have by virtue of Proposition \ref{s2}, \eqref{e7}, and \eqref{ha2} 
that
\be{c6a}
\max\left\{\Var_{[0,T]} \xi\ok, \Var_{[0,T]} \xi\ol\right\} \le 3\Var_{[0,T]} u + 3L \Var_{[0,T]} w =: \hat C.
\ee
Now let us choose $z = \hat x\ok_j \in Z(w\ol_j)$ in \eqref{xxi} for $x\ol_j$ such that 
$|\hat x\ok_j- x\ok_j|=\dist(x\ok_j,Z(w\ol_j))$ and, similarly, $z = \hat x\ol_j \in Z(w\ok_j)$ in \eqref{xxi} for $x\ok_j$ such that $|\hat x\ol_j- x\ol_j|= \dist(x\ol_j,Z(w\ok_j))$, we obtain that
\be{c7}
\begin{aligned}
\scal{\xi\ol_j - \xi\ol_{j-1}, x\ol_j - \hat x\ok_j} + \frac{|\xi\ol_j - \xi\ol_{j-1}|}{2r}\,|x\ol_j - \hat x\ok_j|^2 &\ge 0,\\ 
\scal{\xi\ok_j - \xi\ok_{j-1}, x\ok_j - \hat x\ol_j} + \frac{|\xi\ok_j - \xi\ok_{j-1}|}{2r}\,|x\ok_j - \hat x\ol_j|^2 &\ge 0.
\end{aligned}
\ee 
Put
$$
\Delta := \sup_{j=1, \dots, m} \left(|u\ol_j - u\ok_j| + d_H(Z(w_j\ol),Z(w_j\ok))\right).
$$
We have
$$
\begin{aligned}
\scal{\xi\ol_j - \xi\ol_{j-1}, x\ol_j - \hat x\ok_j} &\le \scal{\xi\ol_j - \xi\ol_{j-1}, x\ol_j - x\ok_j} + |\xi\ol_j - \xi\ol_{j-1}|\,d_H(Z(w_j\ol),Z(w_j\ok)),\\
\scal{\xi\ok_j - \xi\ok_{j-1}, x\ok_j - \hat x\ol_j} &\le \scal{\xi\ok_j - \xi\ok_{j-1}, x\ok_j - x\ol_j} + |\xi\ok_j - \xi\ok_{j-1}|\,d_H(Z(w_j\ol),Z(w_j\ok)),\\
|x\ol_j - \hat x\ok_j|^2 &\le \left(\Delta + |\xi\ol_j - \xi\ok_j|\right)^2 \le (1+\kappa)|\xi\ol_j - \xi\ok_j|^2 + \left(1+\frac1{\kappa}\right) \Delta^2,\\
|\hat x\ol_j - x\ok_j|^2 &\le \left(\Delta + |\xi\ol_j - \xi\ok_j|\right)^2 \le (1+\kappa)|\xi\ol_j - \xi\ok_j|^2 + \left(1+\frac1{\kappa}\right) \Delta^2
\end{aligned}
$$
where $\kappa>0$ is chosen in such a way that
\be{kappa}
r^* = \frac{r}{M(1+\kappa)^2}.
\ee
Furthermore, put
$$
\delta_j = |\xi\ol_j - \xi\ol_{j-1}|+|\xi\ok_j - \xi\ok_{j-1}|
$$
We sum up the inequalities in \eqref{c7} and use the above formulas to obtain
\begin{align}\nonumber
&\scal{(\xi\ol_j - \xi\ok_j) - (\xi\ol_{j-1} - \xi\ok_{j-1}), \xi\ol_j - \xi\ok_j}
 \le \scal{(\xi\ol_j - \xi\ol_{j-1}) - (\xi\ok_j - \xi\ok_{j-1}), u\ol_j - u\ok_j} \\ \nonumber
&\quad + \delta_j\left(\frac{1+\kappa}{2r}|\xi\ok_j {-} \xi\ol_j|^2 + 
 \frac{(1+\kappa)}{2r\kappa}\Delta^2 + d_H(Z(w_j\ol),Z(w_j\ok))\right)\\ \label{c8}
&\le \frac{\delta_j(1+\kappa)}{2r} \left(|\xi\ok_j - \xi\ol_j|^2+ U\right),
\end{align}
where 
$$
U = \frac{1}{\kappa} \Delta^2 + \frac{2r}{1+\kappa}\Delta.
$$Similarly as in \eqref{e12} we have
$$
\scal{(\xi\ol_j - \xi\ok_j) - (\xi\ol_{j-1} - \xi\ok_{j-1}), \xi\ol_j - \xi\ok_j} \ge 
\frac12 \left(|\xi\ol_j - \xi\ok_j|^2 - |\xi\ol_{j-1} - \xi\ok_{j-1}|^2\right).
$$
In order to simplify further the argument, we put
\be{yj}
Y_j := |\xi\ol_j - \xi\ok_j|^2, \quad \lambda_j := \frac{\delta_j(1+\kappa)}{r}.
\ee
Then \eqref{c8} can be written in the form
\be{c9}
Y_j - Y_{j-1} \le \lambda_j (U + Y_j) \ \for j=1, \dots, m.
\ee
By \eqref{c6}--\eqref{c6a} we have
\be{lamj}
\sum_{j=1}^m \lambda_j \le \frac{2C(1+\kappa)}{r},
\ee
where 
$C:=\max\{\bar C, \hat C\}$, and, by \eqref{e7}, \eqref{e5-k}, and \eqref{kappa} we have
\be{c9a}
\delta_j \le \frac{2(4M-1)}{2M-1}r^* = \frac{2r(4M-1)}{M(2M-1)(1+\kappa)^2}.
\ee
Then we get from \eqref{yj} and \eqref{e5g} that
\be{lambda}
\lambda_j \le \frac{2(4M-1)}{M(2M-1)(1+\kappa)} = \frac{1}{1+\kappa}.
\ee
for all $j=1, \dots m$. Then
\begin{equation}\label{c9-1}
 Y_j\prod_{i=1}^j (1-\lambda_i) - Y_{j-1}\prod_{i=1}^{j-1}(1-\lambda_i) 
\le \lambda_j U\prod_{i=1}^{j-1}(1-\lambda_i)  \le \lambda_j U.
\end{equation}
Summing up the above inequality over $j=1, \dots, p$ for an arbitrary $p \le m$ and using the fact that $Y_0 = 0$ 
we obtain
$$
Y_p\prod_{i=1}^p (1-\lambda_i)  \le \sum_{j=1}^p \lambda_j U \le \frac{2C(1+\kappa)}{r} U.
$$
Hence, we obtain for all $p=1, \dots, m$ the estimate
$$
Y_p \le \frac{2C(1+\kappa)}{r} U\prod_{i=1}^p \frac{1}{1-\lambda_i}.
$$
Using \eqref{lamj}, \eqref{lambda}, and the elementary inequality $\log(1+s) \le s$ for $s\ge 0$, we therefore have
$$
\log\left(\prod_{i=1}^m \frac{1}{1-\lambda_i}\right) = \sum_{i=1}^m \log\left(\frac1{1-\lambda_i}\right) 
\le \sum_{i=1}^m \frac{\lambda_i}{1-\lambda_i} \le \frac{1+\kappa}{\kappa}
\sum_{i=1}^m \lambda_i \le \frac{2C(1+\kappa)^2}{\kappa r}.
$$
This yields the final estimate
\be{c10}
\begin{aligned}
&\|\xi\ok - \xi\ol\|^2 = \max\{Y_p;p=1, \dots, m\} \le \tilde C U\\
&\ \le C^*\left(\|u\ok {-} u\ol\|+ \|u\ok {-} u\ol\|^2 + d_H(Z(w\ol),Z(w\ok)) + d_H^2(Z(w\ol),Z(w\ok))\right)
\end{aligned}
\ee
where $\tilde C, C^* > 0$ are suitable constants independent of $k$ and $l$. 
Then $\{\xi\ok\}$ is a Cauchy sequence in the space of right-continuous regulated functions
which has uniformly bounded variation, and we conclude that there exists a function $\xi \in BV_R(0,T;X)$ 
such that $\|\xi\ok - \xi\| \to 0$ as $k \to \infty$.
Then also the functions $x\ok = u\ok - \xi\ok$ converge uniformly to $x = u-\xi$, and $x(t) \in Z(w(t))$ for all 
$t \in [0,T]$ 
and we have
\begin{equation}\label{c10-1}
\int_\sigma^\tau \scal{x^{(k)}(t) - z(t), \dd\xi^{(k)}(t)} + \frac1{2r}\int_\sigma^\tau |x^{(k)}(t) - z(t)|^2 \dd V(\xi^{(k)})(t) \ge 0,
\end{equation}
whenever $0 \le \sigma \le \tau \le T$ and $z \in G(0,T;X)$ with $z(t) \in Z(w\ok(t))$ for all $t\in [0,T]$.
The functions $V\ok(t) := V(\xi\ok)(t)$ are nondecreasing and uniformly bounded by $\max\{\bar{C},\hat{C}\}$ due to \eqref{c6}--\eqref{c6a}; hence, thanks to the uniform convergence of $\xi^{(k)}$, they converge uniformly to a bounded nondecreasing function $\bar V$ and we have that 
$V(\xi)(\tau) - V(\xi)(\sigma)\le \bar V(\tau) - \bar V(\sigma)$ for $0 \le \sigma < \tau \le T$, by the lower semicontinuity of the variation. Therefore passing to the limit as $k \to \infty$ in \eqref{c10-1} yields
\begin{equation}\label{ca2}
\int_\sigma^\tau \scal{x(t) - z(t), \dd\xi(t)} + \frac1{2r}\int_\sigma^\tau |x(t) - z(t)|^2 \dd \bar V(t) \ge 0
\end{equation}
whenever $0 \le \sigma \le \tau \le T$ and $z \in G(0,T;X)$ with $z(t) \in Z(w(t))$ for all $t\in [0,T]$.
By virtue of \cite[Lemma A.9]{Rec11a} we know that the integrals can also be interpreted as Lebesgue integrals over $[\sigma,\tau]$ (see Remark \ref{KS-LS} below) with respect to the differential measures generated by $\xi$ and $\bar{V}$, which we denote $\dd\xi$ and $\dd\bar{V}$, respectively. 
In what follows, given a measure $\mu$, we write $\int_{J}\scal{f,\dd\mu}$ to denote the Lebesgue integral of $f$ with respect to the measure $\mu$ over the interval $J$; and by $L^1(\mu;E)$ we denote the space of $\mu$-integrable $E$-valued functions.
Since $V(\xi)(\tau) - V(\xi)(\sigma)\le \bar V(\tau) - \bar V(\sigma)$ whenever $0 \le \sigma < \tau \le T$, we have that $\dd V(\xi)$ is absolutely continuous with respect to $\dd \bar V$, hence there exists $h \in L^1(\dd\bar V,\mathbb{R})$ such that $\dd V(\xi) = h \dd\bar V$ and $0 \le h(t) \le 1$ for all $t \in [0,T]$. Moreover, $\dd V(\xi)$ is the total variation measure of $\dd\xi$, thus there exists 
$\eta \in L^1(\dd V(\xi);X)$ such that $|\eta(t)| = 1$ for every $t \in [0,T]$ and $\dd\xi = \eta\dd V(\xi)$.
Therefore we infer that $\dd \xi = h\eta\dd \bar V$ and \eqref{ca2} reads
\begin{equation}\label{ca3}
\int_{[\sigma,\tau]} \scal{x(t) - z(t), h(t)\eta(t)} \dd\bar V(t)+ \frac1{2r}\int_{[\sigma,\tau]} |x(t) - z(t)|^2 \dd \bar V(t) \ge 0
\end{equation}
whenever $0 \le \sigma \le \tau \le T$ and $z \in G(0,T;X)$ with $z(t) \in Z(w(t))$ for all $t\in [0,T]$.
Now we divide \eqref{ca3} by $\bar V(\tau) - \bar V(\sigma)$ for every $\dd \bar V$-Lebesgue point $\sigma$ of both integrands $\scal{x(\cdot) - z(\cdot), h(\cdot)\eta(\cdot)}$ and $|x(\cdot) - z(\cdot)|^2$, and for every 
$\tau \in (\sigma,T]$ such that $\bar V(\tau) - \bar V(\sigma) > 0$. Letting $\tau$ go to $\sigma$ we infer that 
\begin{equation}\label{ca4}
\scal{x(t) - z(t), h(t)\eta(t)}+ \frac{1}{2r} |x(t) - z(t)|^2 \ge 0 \qquad \text{for $\dd \bar V$-a.e. $t \in [0,T]$.}
\end{equation} 
If $B_+ \subseteq [0,T]$ is the Borel set where $h(t) > 0$, then from \eqref{ca4} we infer that
\begin{equation}\label{ca5}
\scal{x(t) - z(t), \eta(t)}+ \frac{1}{2rh(t)} |x(t) - z(t)|^2 \ge 0 \qquad \text{for $\dd \bar V$-a.e. $t \in B_+$,}
\end{equation}
which implies that $\eta(t) \in N_{Z(w(t))}(x(t))$ for $\dd \bar V$-a.e. $t \in B_+$ (argue, e.g., as in \eqref{e3f}-\eqref{e3g} recalling that $0 \le h(t) \le 1$). Therefore, as $Z(w(t))$ is 
$r$-prox-regular and $|\eta(t)| = 1$, it follows that 
\begin{equation}\label{ca6}
\scal{x(t) - z(t), \eta(t)}+ \frac{1}{2r} |x(t) - z(t)|^2 \ge 0 \qquad \text{for $\dd \bar V$-a.e. $t \in B_+$}.
\end{equation}
Multiplying \eqref{ca6} by $h(t)$ we therefore obtain that
\begin{equation}\label{ca7}
\scal{x(t) - z(t), h(t)\eta(t)}+ 
\frac{1}{2r} |x(t) - z(t)|^2 h(t)\ge 0 \qquad \text{for $\dd \bar V$-a.e. $t \in [0,T]$},
\end{equation}
hence integrating with respect to $\dd \bar V$ over $[0,T]$ and recalling that $\dd \xi = h\eta\dd \bar V$ and 
$\dd V(\xi) = h \dd \bar V$, we get
\begin{equation}\label{ca8}
\int_{[0,T]}\scal{x(t) - z(t), \dd \xi(t)} + \frac{1}{2r} \int_{[0,T]} |x(t) - z(t)|^2 \dd V(\xi)(t) \ge 0
\end{equation} 
for every $z \in G(0,T;X)$ with $z(t) \in Z(w(t))$ for all $t\in [0,T]$, (where the integrals exist as Lebesgue integrals, therefore also as Kurzweil integrals) i.e. \eqref{s4} holds and the existence part is proved.
Uniqueness follows from Lemma \ref{unique}.
\epf

Theorem \ref{t1} states that the time discrete approximations of $\xi$ defined by the catching up algorithm \eqref{xxi}
converge to the unique solution $\xi$ of \eqref{s4} uniformly.
We now show that in Case (i) of Theorem \ref{t1}, also the input-output relation defined by \eqref{s4} is continuous with respect to the sup-norm. We need two preliminary results.

\begin{lemma}\label{cp}
Assume that $Z_1, Z_2$ are two non-empty closed r-prox-regular sets of $X$ such that 
$d_H(Z_1,Z_2) < \infty$ and assume that $0 <\delta \le r/2$. Then there exists a constant $C > 0$ depending only on $r$ such that for every $y_1, y_2 \in X$, $\dist(y_i,Z_i) \le \delta$ for $i=1,2$, we have 
\begin{equation}\label{cp1}
  |\zeta_1 - \zeta_2|^2 \le C\left(|y_1 - y_2|^2 + d_H^2(Z_1,Z_2) + d_H(Z_1,Z_2)\right)
\end{equation} 
where $\zeta_i \in Z_i$ denotes the unique vector in $Z_i$ such that $|y_i - \zeta_i| = \dist(y_i,Z_i)$ for 
$i=1,2$.
 \end{lemma}
\bpf{Proof}
Using the $r$-prox-regularity of the sets $Z_1$, $Z_2$ we infer that for every $z_1 \in Z_1$ and 
$z_2 \in Z_2$ we have
\begin{align*} 
   |\zeta_1 - \zeta_2|^2
    & =  \scal{y_1 - \zeta_1, z_1 -\zeta_1} +
                         \scal{y_2 - \zeta_2, z_2 - \zeta_2} \\
    & \qquad + \scal{y_1 - \zeta_1, \zeta_2 - z_1} 
                     + \scal{y_2 - \zeta_2,\zeta_1 - z_2} + \scal{y_1 - y_2, \zeta_1 - \zeta_2}  \\
    & \le \frac{1}{2r}|y_1 - \zeta_1| |z_1 -\zeta_1|^2 + \frac{1}{2r}|y_2 - \zeta_2||z_2 - \zeta_2|^2 \notag \\
    & \qquad + |y_1 - \zeta_1||\zeta_2 - z_1| + |y_2 - \zeta_2||\zeta_1 - z_2| + |y_1 - y_2||\zeta_1 - \zeta_2|
    \\
    & \le \frac{\delta}{2r}(|z_1 -\zeta_1|^2 + |z_2 - \zeta_2|^2)  + 
           \delta(|\zeta_2 - z_1| + |\zeta_1 - z_2|)+ |y_1 - y_2||\zeta_1 - \zeta_2|. 
\end{align*}
We now take in the previous inequality $z_1$ and $z_2$ such that 
$\max\{|z_1 - \zeta_2|,|z_2 - \zeta_1|\} \le d_H(Z_1,Z_2)$, and denote $E = |\zeta_1 - \zeta_2|$, 
$D = d_H(Z_1,Z_2)$, $Y = |y_1 - y_2|$. We have $|z_1 -\zeta_1|\le |z_1 -\zeta_2|+|\zeta_2 -\zeta_1|$, 
$|z_2 -\zeta_2|\le |z_2 -\zeta_1|+|\zeta_1 -\zeta_2|$, and from the Young inequality $ab \le a^2/4 + b^2$ for $a,b\in \real$ we get 
\begin{align*}
E^2 &\le \frac{\delta}{r}(E + D)^2 + 2\delta D + YE \le \frac12 (E + D)^2 + r D + YE\\
& \le \frac34 E^2 + (D+Y)^2 + \frac12 D^2 + rD,
\end{align*}
and \eqref{cp1} follows.
\epf

\begin{corollary}\label{cp3}
Let $\{Z(w); w \in A\}$ be a family of $r$-prox-regular sets satisfying \eqref{ha1}. Assume that 
$\delta > 0$, $y \in G_R(0,T;X)$, $w \in G_R(0,T;X)$ are given such that $w(t) \in A$ and 
$\dist(y(t),Z(w(t))) \le r/2$ for every $t\in[0,T]$. Let $\zeta : [0,T] \to X$ be defined in such a way that
$\zeta(t)$ is the only vector in $Z(w(t))$ such that
\begin{equation}
  \scal{y(t) - \zeta(t), \zeta(t) - z} + \frac{|y(t) - \zeta(t)|}{2r}|\zeta(t) - z|^2  \ge 0 \qquad \forall z \in Z(w(t)).
\end{equation}
Then $\zeta \in G_R(0,T;X)$.
\end{corollary}

\bpf{Proof}
Let us fix $t_0 \in (0,T]$. For each $\ve>0$ we use \eqref{ha1} to find $t_1<t_0$ such that
\be{left}
|y(\tau) - y(t_0-)| + d_H(Z(w(\tau)),Z(w(t_0-))) < \ve \ \for \tau \in (t_1, t_0).
\ee
Thanks to Lemma \ref{cp} there exists a constant $C > 0$ such that 
\[
  |\zeta(t) - \zeta(s)|^2 \le C(|y(t) - y(s)|^2 + d_H(Z(w(t)),Z(w(s))) + d_H^2(Z(w(t)),Z(w(s)))
\]
whenever $t_1 < t < s < t_0$, so the existence of $\zeta(t_0-) \in X$ follows from \eqref{left}. The argument for the right limits is analogous.
\epf

\begin{theorem}\label{t2}
Let $\{Z(w); w \in A\}$ be a family of $r$-prox-regular sets satisfying \eqref{ha1} and having 
uniformly non-empty interior with $\rho, R$ as in Definition \ref{d2},
and consider the subset $D\subset G_R(0,T;X)\times G_R(0,T;W) \times X$ defined by
$$D=\{(u,w,x_0): w(t)\in A \mbox{ for }t\in[0,T],\,x_0\in Z(w(0)),\mbox{and \eqref{e5g} holds}\}$$
Then the mapping $\mathcal{R}: (u,w,x_0)\in D \mapsto \xi\in G_R(0,T;X)$ 
which with given $u \in G_R(0,T;X)$, $w\in G_R(0,T;W)$ and an initial condition $x_0 \in Z(w(0))$
associates the solution $\xi \in G_R(0,T;X)$ is continuous with respect to the norm $\|\cdot\|$. 
\end{theorem}

\bpf{Proof}
Consider sequences $\{u_n\}$ in $G_R(0,T;X)$ and $\{w_n\}$ in $G_R(0,T;W)$, $w_n(t) \in A$ for all $t \in [0,T]$ which converge 
 uniformly to $u$ and $w$, respectively, and initial conditions $x^0_n\in Z(w_n(0))$ which converge to $x_0$. We proceed as in the proof of Theorem \ref{t1} and for $\rho$ as in Proposition \ref{p1} we find a division $0 = \hat t_0 < \hat t_1 < \dots < \hat t_N = T$ such that
\be{c14}
d_H(Z(w(t)), Z(w(\hat t_{i-1}))) + |u(t) - u(\hat t_{i-1})| \le \frac\rho{9} \ \for t \in t \in [\hat t_{i-1}, \hat t_i), \ i=1, \dots N.
\ee
We further find $n_0 \in \nat$ such that for $n \ge n_0$ we have
\be{c15}
\sup_{t\in [0,T]}d_H(Z(w(t)), Z(w_n(t))) +\|u - u_n\| < \frac\rho{9}.
\ee
For each $n \ge n_0$ we construct a piecewise constant approximations $u\ok_n$ of $u_n$ and $w\ok_n$ of $w_n$ as in \eqref{c4}-\eqref{c4w} and find $k_n\in \nat$ such that 
$$
\sup_{t\in [0,T]}d_H(Z(w_n(t)), Z(w\ok_n(t))) +\|u_n - u\ok_n\| < \frac\rho{3} \ \for k\ge k_n.
$$
By the triangle inequality we have for all $n\ge n_0$ and $k\ge k_n$ that
\be{c16}
d_H(Z(w\ok_n(t)), Z(w\ok_n(\hat t_{i-1}))) + |u\ok_n(t) - u\ok_n(\hat t_{i-1})| \le \rho \ \for t \in [\hat t_{i-1}, \hat t_i),, \ i=1, \dots N.
\ee
We now proceed as in the proof of Theorem \ref{t1} and obtain for the piecewise constant solutions $\xi\ok_n$, as a counterpart of \eqref{c6}, the estimate
\be{c17}
\sup_{n\ge n_0, k\ge k_n} \Var_{[0,T]} \xi\ok_n \le \bar C
\ee
with a constant $\bar C$ depending only on $u$ and the geometry of the sets $Z(w)$. We already know that for every $n \ge n_0$, $\xi\ok_n$ converge uniformly as $k \to \infty$ to the unique solution $\xi_n$ of the variational inequality
\be{c18}
\int_0^t \scal{u_n(\tau) - \xi_n(\tau) - z(\tau), \dd\xi_n(\tau)} + 
\frac1{2r}\int_0^t |u_n(\tau) - \xi_n(\tau) - z(\tau)|^2 \dd V(\xi_n)(\tau) \ge 0
\ee
for every $t \in [0,T]$ and every $z \in G(0,t;X)$, $z(\tau) \in Z(w_n(\tau))$ for $\tau\in [0,t]$, where we have also used Lemma \ref{kl1}. Moreover, we have
\be{c17a}
\sup_{n\ge n_0} \Var_{[0,T]} \xi_n \le \bar C
\ee
by virtue of \eqref{c17}. On the other hand, we have
\be{c19}
\int_0^t \scal{u(\tau) - \xi(\tau) - z(\tau), \dd\xi(\tau)} + 
\frac1{2r}\int_0^t |u(\tau) - \xi(\tau) - z(\tau)|^2 \dd V(\xi)(\tau) \ge 0
\ee
for every $t \in [0,T]$ and every $z \in G(0,t;X )$, $z(\tau) \in Z(w(\tau))$ for $\tau\in [0,t]$.
Let us observe now that for every $t \in [0,T]$ 
we have  $u(t)-\xi(t) \in Z(w(t))$, and we find $\zeta_n(t) \in Z(w_n(t))$ such that 
$|u(t) - \xi(t) - \zeta_n(t)| = \dist(u(t) - \xi(t), Z(w_n(t)))$. Similarly, we find $\bar{\zeta}_n(t) \in Z(w(t))$ such that 
$|u_n(t)-\xi_n(t) - \bar{\zeta}_n(t)| = \dist(u_n(t) - \xi_n(t), Z(w(t)))$. 
Let us set for simplicity $\Delta_n = \sup_{\tau} d_H(Z(w(\tau)), Z(w_n(\tau)))+ \|\bar u_n\|$, 
$\bar u_n = u-u_n$ and $\bar{\xi}_n=\xi-\xi_n$.  From Corollary \ref{cp3} it follows that $\bar{\zeta}_n$ and $\zeta_n$ are regulated for every $n$ sufficiently large,
hence, putting $z(\tau) = \bar{\zeta}_n(\tau)$ in \eqref{c19}, we obtain
\begin{align}\label{c20}
  & -\int_0^t \scal{u(\tau) - \xi(\tau) - \bar{\zeta}_n(\tau), \dd \xi(\tau)} \notag \\
   & \qquad \le \frac{1}{2r} \int_0^t |u(\tau) - \xi(\tau) - \bar{\zeta}_n(\tau)|^2 \dd V(\xi)(\tau) \notag \\
   & \qquad \le \frac{1}{2r} \int_0^t  (|\bar u_n(\tau)| + |\bar \xi_n(\tau)| + |u_n(\tau) - \xi_n(\tau) -         
      \bar{\zeta}_n(\tau)|)^2 \dd V(\xi) (\tau) \notag \\
   & \qquad \le \frac{1}{2r} \int_0^t (2|\bar \xi_n(\tau)|^2 + 4|\bar u_n(\tau)|^2 + 4d_H^2(Z(w_n(\tau)), Z(w(\tau))))\dd V(\xi)(\tau)
   \notag \\
   & \qquad \le \frac{1}{r} \int_0^t (|\bar \xi_n(\tau)|^2+ 2\Delta_n^2)\dd V(\xi)(\tau)
\end{align}
Putting $z(\tau) = \zeta_n(\tau)$ in \eqref{c18} and using the same argument as in \eqref{c20} yields
\begin{equation}\label{c20-1}
   -\int_0^t \scal{u_n(\tau) - \xi_n(\tau) - \zeta_n(\tau), \dd \xi_n(\tau)} 
    \le \frac{1}{r} \int_0^t (|\bar \xi_n(\tau)|^2+ 2\Delta_n^2)\dd V(\xi_n)(\tau).
\end{equation}
On the other hand we have that
\begin{align}\label{c20-2}
 & \int_0^t \scal{\bar \xi_n(\tau), \dd\bar \xi_n(\tau)} \notag \\
 & \qquad =  \int_0^t \scal{ \bar \xi_n(\tau) - \bar u_n(\tau), \dd \bar \xi_n(\tau)} 
                    +\int_0^t \scal{\bar u_n(\tau), \dd \bar \xi_n(\tau)} \notag \\
 & \qquad =  \int_0^t \scal{ \bar \xi_n(\tau) - \bar u_n(\tau), \dd \xi(\tau)} 
                     - \int_0^t \scal{ \bar \xi_n(\tau) - \bar u_n(\tau), \dd \xi_n(\tau)} 
                    + \int_0^t \scal{\bar u_n(\tau), \dd \bar \xi_n(\tau)} \notag \\                   
 & \qquad = \int_0^t \scal{\xi(\tau) - u(\tau) + \bar{\zeta}_n(\tau), \dd \xi(\tau)} 
                    + \int_0^t \scal{u_n(\tau) - \xi_n(\tau) - \bar{\zeta}_n(\tau), \dd \xi(\tau)} \notag \\
 & \qquad\qquad - \int_0^t \scal{\xi(\tau) - u(\tau) + \zeta_n(\tau), \dd \xi_n(\tau)} 
                    - \int_0^t \scal{u_n(\tau) - \xi_n(\tau) - \zeta_n(\tau), \dd \xi_n(\tau)} \notag \\ 
 & \qquad\qquad + \int_0^t \scal{\bar u_n(\tau), \dd \bar \xi_n(\tau)}. 
\end{align}
Therefore from \eqref{c20}, \eqref{c20-1}, and \eqref{c20-2} we infer that
\begin{align}
 & \int_0^t \scal{\bar \xi_n(\tau), \dd\bar \xi_n(\tau)} \notag \\
  & \qquad \le  \int_0^t |u_n(\tau) - \xi_n(\tau) - \bar{\zeta}_n(\tau)| \dd V(\xi)(\tau) 
   + \int_0^t |u(\tau) - \xi(\tau) - \zeta_n(\tau)| \dd V(\xi_n)(\tau)    \notag \\
   & \qquad\qquad  +    \frac{1}{r} \int_0^t (|\bar \xi_n(\tau)|^2 + 2\Delta_n^2 )\dd(V(\xi) + V(\xi_n))(\tau)    
   + \int_0^t |\bar u_n(\tau)| \dd V(\bar \xi_n)(\tau) \notag \\
 &  \qquad\le \int_0^t  \Delta_n  \dd(V(\xi) + V(\xi_n))(\tau)   +   
  \frac{1}{r} \int_0^t (|\bar \xi_n(\tau)|^2 + 2\Delta_n^2 )\dd(V(\xi) + V(\xi_n))(\tau),              
\end{align}
and using Corollary \ref{3c17} we obtain that there exists a suitable constant $C>0$ independent of $n$ such that
\be{c21}
|\bar \xi_n(t)|^2 \le |\bar \xi_n(0)|^2 + 
\int_0^t \left(\Delta_n + \Delta_n^2 + |\bar \xi_n(\tau)|^2\right) \dd g_n(\tau)
\ee
with $g_n(\tau) = C(V(\xi_n)(\tau)+V(\xi)(\tau))$. This is an inequality of the form \eqref{g1z} with $z(t) = \Delta_n + \Delta_n^2 + |\bar \xi_n(t)|^2$. The functions $g_n$ are bounded above independently of $n$ as a consequence of \eqref{c17a}. From Lemma \ref{gl2} it thus follows that there exists a constant $C>0$ independent of $n$ such that
$$
|\bar \xi_n(t)|^2 \le C\left(\Delta_n + \Delta_n^2 + |\bar \xi_n(0)|^2\right),
$$ 
for every $t \in [0,T]$, and the assertion follows easily.
\epf


\section{More regular inputs}\label{con}

In this section we consider the cases when the inputs are continuous or absolutely continuous and prove that so is the solution of Problem \ref{PB} under appropriate assumptions. Before proceeding with the continuous case we first prove a local result.

\begin{lemma}\label{al1x}
Let $\{Z(w); w\in A\} \subset X$ be a family of $r$-prox-regular sets satisfying \eqref{ha1},
and let $u \in G_R(0,T; X)$, $w \in G_R(0,T; W)$, and $d\in (0,r/2]$ be given, and for 
$0 \le \sigma < \tau \le T$ such that for all $t \in [\sigma,\tau]$ we have $d_H(Z(w(t)), Z(w(\sigma)) \le d$
put
$$
\begin{aligned}
U(\sigma,\tau) &= \sup_{t \in [\sigma, \tau]} |u(t) - u(\sigma)| + 
                             \sup_{t \in [\sigma, \tau]} |u(t) - u(\sigma)|^2\\
&\quad +\sup_{t \in [\sigma, \tau]} d_H(Z(w(t)), Z(w(\sigma))) + 
               \sup_{t \in [\sigma, \tau]} d_H^2(Z(w(t)), Z(w(\sigma))).
\end{aligned}
$$
Let $\xi \in BV_R(0,T;X)$ be a function satisfying \eqref{s4}--\eqref{s5}.
Then there exists a constant $C^* > 0$ such that for all $0 \le \sigma < s \le \tau \le T$ we have
\be{a3x}
|\xi(s) - \xi(\sigma)|^2 \le C^* U(\sigma,\tau) \Var_{[\sigma,\tau]}(\xi).
\ee
\end{lemma}

\bpf{Proof}
Let $x = u - \xi$. By Lemma \ref{kl1} for arbitrary $0 \le \sigma < s \le \tau \le T$, where we choose 
$\tilde z(t) = \hat x(\sigma,t)\in Z(w(t))$ for $t\in [\sigma,\tau]$ in \eqref{k5} such that 
$|\hat x(\sigma,t) - x(\sigma)| = \dist(x(\sigma), Z(w(t)))$. 
By virtue of Corollary \ref{cp3}, this is an admissible choice, and we obtain
\be{a2p}
\int_\sigma^s \scal{u(t) - \xi(t) - \hat x(\sigma,t), \dd \xi(t)} + \frac{1}{2r}\int_\sigma^s |u(t) - \xi(t) - \hat x(\sigma,t)|^2\dd V(\xi)(t) \ge 0,
\ee
hence,
\be{a2}
\begin{aligned}
&\int_\sigma^s \scal{u(t) - \xi(t) - u(\sigma) + \xi(\sigma), \dd \xi(t)}- \int_\sigma^s \scal{\hat x(\sigma,t) - x(\sigma), \dd \xi(t)}\\
&\quad + \frac{1}{2r}\int_\sigma^s |u(t) - \xi(t) - u(\sigma) + \xi(\sigma) - \hat x(\sigma,t) + x(\sigma)|^2\dd V(\xi)(t)  \ge 0.
\end{aligned}
\ee
{}From this inequality and from Corollary \ref{3c17} we infer that
\begin{align*}
\frac12|\xi(s) - \xi(\sigma)|^2 & \le \int_\sigma^s \scal{\xi(t) {-} \xi(\sigma), \dd \xi(t)} \\ 
&\le \int_\sigma^s \scal{u(t) {-} u(\sigma), \dd \xi(t)}- \int_\sigma^s \scal{\hat x(\sigma,t) - x(\sigma), \dd \xi(t)}\\
&\qquad+ \frac{3}{2r}\int_\sigma^s \left(|u(t) {-} u(\sigma)|^2 + |\xi(t) {-} \xi(\sigma)|^2 + |x(\sigma) - \hat x(\sigma,t)|^2\right)\dd V(\xi)(t)\\
&\le C \int_\sigma^s \left(U(\sigma,\tau) + |\xi(t) - \xi(\sigma)|^2\right)\dd V(\xi)(t)
\end{align*}
with $C = \max\{1, 3/(2r)\}$.
We now use Lemma \ref{gl2}  with $\gamma= 2CU(\sigma,\tau)\Var_{[\sigma,\tau]}\xi$, 
$z(t) = |\xi(t+\sigma) - \xi(\sigma)|^2$ and $g(t) = 2C V(\xi)(t+\sigma)$ for $t \in [0, \tau - \sigma]$. 
Lemma \ref{gl2} yields $z(t)\le \gamma y(t)$ for $t\in [0, \tau - \sigma]$. Since $y$ is bounded by virtue of \eqref{g9a-bis} of Lemma \ref{gl1}
and the bounded variation of $\xi$, we obtain
\be{a3}
|\xi(s) - \xi(\sigma)|^2 \le 2C\expe^{2C\Var_{[0,T]}(\xi)}U(\sigma,\tau) \Var_{[\sigma,\tau]}(\xi),
\ee
which we wanted to prove.
\epf

As a first corollary of Lemma \ref{al1x} we prove that in Case (i) of Theorem \ref{t1}, the output is continuous if the inputs are continuous. We thus provide an independent proof of the result obtained in \cite[Theorem 4.2]{cmm} reformulated in terms of the Kurzweil integral.

\begin{corollary} \label{cx1}
Let $\{Z(w); w\in A\} \subset X$ be a family of $r$-prox-regular sets satisfying \eqref{ha1}  and having 
uniformly non-empty interior with $\rho, R$ as in Definition \ref{d2}, and let $u \in C([0,T];X)$, $w\in C([0,T];A)$, and $x_0 \in Z(w(0))$ be given.
Then there exists a unique $\xi \in BV(0,T;X) \cap C([0,T];X)$ and a unique $x \in C([0,T];X)$, 
$x(t) \in Z(w(t))$, $\xi(t) + x(t) = u(t)$ for all $t\in [0,T]$, such that \eqref{s4} holds and $x(0) = x_0$.
\end{corollary}

\bpf{Proof}
The existence and uniqueness of the solution $\xi$ follows from Case (i) of Theorem \ref{t1}, and the continuity of $\xi$ is an immediate consequence of Lemma \ref{al1x}.
\epf

Finally we turn our attention to the absolutely continuous case. The result reads as follows.

\begin{corollary}\label{al1}
Let $\{Z(w); w\in A\} \subset X$ be a family of $r$-prox-regular sets satisfying \eqref{ha2}, 
and let $u \in W^{1,1}(0,T; X)$, $w \in W^{1,1}(0,T; W)$ be arbitrarily chosen. Then the solution $\xi$ to \eqref{s4} with any initial condition $x_0 \in Z(w(0))$ belongs to $W^{1,1}(0,T; X)$ and the variational inequality
\be{a1}
\scal{x(t) - z, \dot\xi(t)} + \frac{|\dot\xi(t)|}{2r}|x(t) - z|^2 \ge 0, \quad x(t) + \xi(t) = u(t),
\ee
holds for a.\,e. $t \in (0,T)$ and all $z \in Z(w(t))$.
\end{corollary}

\bpf{Proof}
Let us start by observing that by virtue of Case (ii) of Theorem \ref{t1} there exists a unique 
$\xi \in BV_R(0,T;X)$ satisfying \eqref{s4}--\eqref{s5} with $x(t) = u(t) - \xi(t)$ for every $t \in [0,T]$.
As in Lemma \ref{al1x} we consider $0 \le \sigma < s \le \tau \le T$ such that 
$\Var_{[\sigma,\tau]}(w)<r/L$. Using \eqref{ha2} we obtain
$$
U(\sigma,\tau) \le (1+\|u\|)\Var_{[\sigma,\tau]}(u)
+ C(1+\|w\|_W)\Var_{[\sigma,\tau]}(w)
$$
with a constant $C$ depending only on $L$. By \eqref{a3x} there exists a constant $\hat C >0$ such that for all $0\le \sigma < \tau \le T$ we have
\be{a4}
|\xi(\tau) - \xi(\sigma)| \le \hat C \sqrt{\Var_{[\sigma,\tau]}(\xi)\,\left(\Var_{[\sigma,\tau]}(u)+\Var_{[\sigma,\tau]}(w)\right)}.
\ee
For any $0\le a < b \le T$ and any sufficiently fine division $a=t_0<t_1<\dots < t_m = b$ such that 
$\Var_{[t_{j-1}, t_j]}(w)\le d < r$ for all $j=1, \dots, m$ we thus have by the Cauchy-Schwarz inequality 
\begin{align*}
\sum_{j=1}^m |\xi(t_j) - \xi(t_{j-1})| &\le
\hat C \sum_{j=1}^m\sqrt{\Var_{[t_{j-1}, t_j]}(\xi)\left(\Var_{[t_{j-1}, t_j]}(u) + \Var_{[t_{j-1}, t_j]}(w)\right)}\\
& \le \hat C \sqrt{\sum_{j=1}^m\Var_{[t_{j-1}, t_j]}(\xi)}\sqrt{\sum_{j=1}^m\left(\Var_{[t_{j-1}, t_j]}(u) + \Var_{[t_{j-1}, t_j]}(w)\right)},
\end{align*}
hence,
$$
\Var_{[a,b]}(\xi) \le \hat C \sqrt{\Var_{[a,b]}(\xi)\,\left(\Var_{[a,b]}(u)+\Var_{[a,b]}(w)\right)}
$$
which implies
\be{a5}
\Var_{[a,b]}(\xi) \le \hat C^2 \left(\Var_{[a,b]}(u) + \Var_{[a,b]}(w)\right)= \hat C^2 \int_a^b \left(|\dot u(t)| + |\dot w(t)|_W\right)\dd t.
\ee
Since \eqref{a5} holds for all $a<b$, we conclude that $\xi$ is absolutely continuous and $|\dot \xi(t)| \le \hat C^2 (|\dot u(t)|+|\dot{w}(t)|)$ almost everywhere. Inequality \eqref{a1} then follows from the general theory of the Kurzweil integral and from Lemma \ref{kl1}.
\epf


\appendix
\section{Appendix}\label{appe}

In this section we recall some basic facts about the Kurzweil integral which are needed throughout the paper.  A good account on such a theory, though restricted to integration of real-valued functions, is the monograph \cite{MST}. The results here are stated for functions with values in the space $X$ endowed with a scalar product $\scal{\cdot,\cdot}$.
Analogous statements in $X$ are proved for the Young integral in \cite{KreLau02} and \cite{bk03}. Note that under the hypotheses of Theorem \ref{Ap5} below, the Kurzweil and the Young integral coincide, see \cite{negli}.

Let $[a,b]$ be a nondegenerate interval of $\real$ and let $\Gamma(a,b)$
be the set of strictly positive functions on $[a,b]$: any $\delta \in \Gamma(a,b)$ is called \emph{gauge} in the framework of the Kurzweil integration. A \emph{partition associated with a division 
$a = t_0 < t_1 < \cdots < t_m = b$} is a set of the form
\begin{equation}\label{Ap1}
  D = \{(\tau_j, [t_{j-1},t_j]);\ \tau_j \in [t_{j-1},t_j],\ j = 1,\ldots,m\},
\end{equation}
and if $\delta \in \Gamma(a,b)$ we say that $D$ is $\delta$-fine if 
\begin{equation}\label{Ap2}
\begin{cases}
   [t_{j-1},t_j] \subset (\tau_j - \delta(\tau_j), \tau_j + \delta(\tau_j)) & \text{for $j = 1,\ldots, m$},\\
   t_{j-1} < \tau_j &   \text{for $j = 2,\ldots, m$}, \\
   \tau_j < t_j &  \text{for $j = 1,\ldots, m-1$}.  
\end{cases}
\end{equation}
It can be proved that the set $\mathcal{F}_\delta(a,b)$ of $\delta$-fine partitions is non-empty, so that
for $f : [a,b] \to X$ and $g : [a,b] \to X$ given, we can define the \emph{Kurzweil integral sum}
\begin{equation}\label{Ap3}
  K_D(f,g) := \sum_{j=1}^m \scal{f(\tau_j), g(t_j) - g(t_{j-1})}.
\end{equation}
We say that $J \in \real$ is the \emph{Kurzweil integral over $[a,b]$ of $f$ with respect to $g$} if for every $\ve > 0$ there exists a $\delta \in \Gamma(a,b)$ such that for every $D \in \mathcal{F}_\delta(a,b)$ we have $|J - K_D(f,g)| < \ve$. In this case we write
\begin{equation}\label{Ap4}
  J = \int_a^b \scal{f(t),\dd g(t)},
\end{equation}
and if $X$ is the real line $\real$ we consequently write $J = \int_a^b f(t) \dd g(t)$. It is easily seen that 
if $J$ exists, then it is unique.

Here is a sufficient condition for the existence of the Kurzweil integral which is sufficient to our purposes
and is tacitly used in the paper.

\begin{theorem}\label{Ap5}
If $f \in G(a,b;X)$ and $g \in BV(0,T;X)$ then $\int_a^b \scal{f(t),\dd g(t)}$ exists. Moreover the 
function $(f,g) \mapsto \int_a^b \scal{f(t),\dd g(t)}$ is bilinear on $G(a,b;X) \times BV(a,b;X)$
\end{theorem}

\begin{remark}\label{KS-LS} It is worth mentioning that for $f \in G(a,b;X)$ and $g \in BV(0,T;X)$, the Kurzweil integral $\int_a^b \scal{f(t),\dd g(t)}$ coincides with the Young integral as well as with the Lebesgue integral with respect to the vector measure generated by $g$. This fact can be easily shown for step functions and extended to regulated functions via convergence theorems. See \cite{KreLau02, Rec11a} for results on variational inequalities involving these integrals. 
\end{remark}

\begin{theorem}\label{add-int}
If the integral $\int_a^b\scal{f(t), \dd g(t)}$ exists, then $\int_c^d\scal{f(t), \dd g(t)}$ exists for every subinterval $[c,d]\in[a,b]$. In particular, for every $c\in(a,b)$
$$\int_a^b\scal{f(t), \dd g(t)}=\int_a^c\scal{f(t), \dd g(t)}+\int_c^b\scal{f(t), \dd g(t)}.$$
\end{theorem}

\begin{theorem}\label{cases}
For any function $f:[a,b]\to X$ and $v\in X$ we have:
\begin{itemize}
\item[$(i)$]  $\dis \int_a^b\scal{f(t),\dd\big(v\chi_{[a,\tau)}\big)(t)}=-\scal{f(\tau),v}$ for $\tau\in(a,b]$,
\item[$(ii)$]  $\dis \int_a^b\scal{f(t),\dd\big(v\chi_{\{ \tau\}}\big)(t)}=\begin{cases}
0 & \mbox{if \ }\tau\in(a,b),
\\
-\scal{f(a),v} & \mbox{if \ }\tau = a,
\\
\scal{f(b),v} & \mbox{if \ } \tau=b.
\end{cases}$
\end{itemize}
\end{theorem}

\begin{theorem}\label{hake}
Let $f:[a,b]\to X$ and $g\in G_R(a,b;X)$ be such that $\int_a^b\scal{f(t), \dd g(t)}$ exists. Then for $c\in(a,b]$ we have
$$\int_a^c\scal{f(t), \dd g(t)}=\lim_{s\to c-}\int_a^s\scal{f(t), \dd g(t)}+\scal{f(c), g(c)-g(c-)}.$$
\end{theorem}

We also need the following convergence result.
\begin{theorem}\label{Ap9}
Assume that $f, f_n \in G(a,b;X)$ and $g, g_n \in BV(a,b;X)$ for every $n \in \nat$. If $\|f-f_n\| \to 0$ and 
$\|g-g_n\| \to 0$ as $n \to \infty$, and if $\sup_n \Var_{[a,b]} g_n < \infty$, then
\begin{equation}\label{a10}
  \lim_{n \to \infty} \int_a^b \scal{f_n(t),g_n(t)} = \int_a^b \scal{f(t),g(t)}.
\end{equation}  
\end{theorem}

For scalar-valued functions, we have the following easy comparison lemma.

\begin{lemma}\label{kl9}
Let $f_1, f_2 \in G(a,b;\real)$ be such that $f_1(t) \le f_2(t)$ for all $t \in [a,b]$, and let $g:[a,b] \to \real$ be nondecreasing. Then
$$
\int_a^b f_1(t) \dd g(t) \le \int_a^b f_2(t) \dd g(t).
$$
\end{lemma}

Integration by parts in the Kurzweil theory involves additional jump terms and the result reads as follows.

\begin{theorem}\label{3t16}
For every $f,g \in BV(a,b;X)$ we have
\begin{gather}
\int_a^b \scal{f(t), \dd g(t)} + \int_a^b \scal{g(t), \dd f(t)} =
\scal{f(b), g(b)} - \scal{f(a), g(a)} \nonumber\\
+ \sum_{t\in [a,b]}\Big(\scal{f(t) - f(t-), g(t) - g(t-)}
- \scal{f(t+) - f(t), g(t+) - g(t)}\Big) .\label{3e28}
\end{gather}
\end{theorem}

\begin{corollary}\label{3c17}
For every $g \in BV_R(a,b;X)$ we have
\be{3e29}
\int_a^b \scal{g(t), \dd g(t)} = \frac 12
|g(b)|^2 - \frac12|g(a)|^2 + \frac 12
\sum_{t\in [a,b]}|g(t) - g(t-)|^2 .
\ee
\end{corollary}

Of course, the number of summands in \eqref{3e28} and \eqref{3e29} is at most countable and the sums are finite.

A Gronwall-type argument exists in the Kurzweil theory, too, but it is less elementary than for the Lebesgue integral. Herein we present and prove a simplified version of such result which is sufficient for our purposes.

\begin{lemma}\label{gl1}
Let $g:[0,T] \to \real$ be a right-continuous nondecreasing function such that for all $t \in (0,T]$ we have
\begin{equation}\label{g-12}
g(t) - g(t-) \le \frac12.
\end{equation}
Then the Kurzweil integral equation
\be{g1}
y(t) = 1 + \int_0^t y(\tau)\dd g(\tau) \ \ \forall t\in [0,T]
\ee
has a unique nondecreasing right-continuous solution $y: [0,T] \to [1,\infty)$ and
\be{g9a-bis}
y(t) \le \,\expe^{2 (g(T)-g(0))} \qquad \forall t \in [0,T].
\ee
\end{lemma}

\bpf{Proof}
Assume first that $g$ is a step function of the form
\be{g2}
g(t) = \sum_{j=1}^{m} g_{j-1} \chi_{[t_{j-1}, t_j)}(t) + g_m \chi_{\{t_m\}}(t)
\ee
with $g_0 \le g_1 \le ... \le g_m$, $g_j - g_{j-1} \le 1/2$ for $j=1, \dots, m$. Then $y$ is a solution of \eqref{g1} if and only if for all $j=1, \dots, m$ and $t \in [t_{j-1}, t_j)$ we have $y(0) = 1$ and
\be{g3}
y(t) = 1 + \sum_{i=1}^{j-1} y(t_i) (g_i - g_{i-1}), \quad y(t_m) = 1 + \sum_{i=1}^{m} y(t_i) (g_i - g_{i-1}).
\ee
In other words, $y$ is a step function of the form
\be{g4}
y(t) = \sum_{j=1}^{m} y_{j-1} \chi_{[t_{j-1}, t_j)}(t) + y_m \chi_{\{t_m\}}(t)
\ee
with
\be{g5}
y_j = 1 + \sum_{i=1}^{j} y_i (g_i - g_{i-1}) \quad \for j = 0,1,\dots, m.
\ee
{}From \eqref{g5} it follows for all $j\ge 1$ that $y_j - y_{j-1} =  y_j (g_j - g_{j-1})$, and we easily conclude by induction that
\be{g6}
y_j = \prod_{i=1}^j \frac{1}{1- g_i + g_{i-1}} \ge 1 \ \for j=1, \dots, m, \quad y_0 = 1.
\ee
Consider now a sequence $\{g\ok; k\in \nat\}$ of nondecreasing  step functions of the form \eqref{g2} which converges uniformly to $g$ as $k\to \infty$ and such that $g\ok(t) - g\ok(t-) \le 1/2$ for all $k \in \nat$, and $g\ok(0) = g(0)$, $g\ok(T) = g(T)$.
We prove that the associated sequence $\{y\ok\}$ of solutions to the equation
\be{g1k}
y\ok(t) = 1 + \int_0^t y\ok(\tau)\dd g\ok(\tau) \ \forall\ t\in [0,T]
\ee
is a Cauchy sequence in $G_R(0,T;\real)$
and the solution $y$ of \eqref{g1} is obtained by passing to the limit as $k \to \infty$ in \eqref{g1k}.

To this end, we find $k_0 \in \nat$ such that
\be{k0}
k \ge k_0 \ \Longrightarrow \ \|g\ok - g\| < \frac{1}{16}.
\ee
Let $k,l \ge k_0$ be fixed, and let $g\ok, g\ol$ be of the form
\be{g2kl}
\begin{aligned}
g\ok(t) &= \sum_{j=0}^{m} g\ok_{j-1} \chi_{[t_{j-1}, t_j)}(t) + g\ok_m \chi_{\{t_m\}}(t),\\
g\ol(t) &= \sum_{j=0}^{m} g\ol_{j-1} \chi_{[t_{j-1}, t_j)}(t) + g\ol_m \chi_{\{t_m\}}(t).
\end{aligned}
\ee
The corresponding solutions $y\ok, y\ol$ are
\be{gy2kl}
\begin{aligned}
y\ok(t) &= \sum_{j=0}^{m} y\ok_{j-1} \chi_{[t_{j-1}, t_j)}(t) + y\ok_m \chi_{\{t_m\}}(t),\\
y\ol(t) &= \sum_{j=0}^{m} y\ol_{j-1} \chi_{[t_{j-1}, t_j)}(t) + y\ol_m \chi_{\{t_m\}}(t),
\end{aligned}
\ee
and we have
\be{g6kl}
\frac{y\ol_j}{y\ok_j} = \prod_{i=1}^j \frac{1- g\ok_i + g\ok_{i-1}}{1- g\ol_i + g\ol_{i-1}} \quad \for j=1, \dots, m.
\ee
To make the formulas short, we denote $\Delta\ok_i = g\ok_i - g\ok_{i-1}$, $\Delta\ol_i = g\ol_i - g\ol_{i-1}$, and rewrite \eqref{g6kl} for $j=1, \dots, m$ as
\begin{align}\nonumber
\log \left(\frac{y\ol_j}{y\ok_j}\right) &= \sum_{i=1}^j \log\left(\frac{1- \Delta\ok_i}{1- \Delta\ol_i}\right)
 = \sum_{i=1}^j \left((\Delta\ol_i {-} \Delta\ok_i) + \log\left(\frac{1- \Delta\ok_i}{1- \Delta\ol_i}\right)-(\Delta\ol_i {-} \Delta\ok_i)\right)\\ \label{g7}
&= g\ol_j - g\ok_j + \sum_{i=1}^j \left(\log\left(\frac{1- \Delta\ok_i}{1- \Delta\ol_i}\right)-(\Delta\ol_i - \Delta\ok_i)\right),
\end{align}
where
$$
\begin{aligned}
&\log\left(\frac{1- \Delta\ok_i}{1- \Delta\ol_i}\right)-(\Delta\ol_i - \Delta\ok_i)
= \log\left(1 + \frac{\Delta\ol_i- \Delta\ok_i}{1- \Delta\ol_i}\right)-(\Delta\ol_i - \Delta\ok_i)\\
&\qquad = \log\left(1 + \frac{\Delta\ol_i- \Delta\ok_i}{1- \Delta\ol_i}\right)-\frac{\Delta\ol_i - \Delta\ok_i}{1- \Delta\ol_i} + \frac{(\Delta\ol_i - \Delta\ok_i)\Delta\ol_i}{1- \Delta\ol_i}.
\end{aligned}
$$
We have $1- \Delta\ol_i \ge 1/2$ and, by \eqref{k0}, $|\Delta\ol_i- \Delta\ok_i|\le 1/4$.
Using the formula $|\log(1+s) - s| \le s^2$ for every $|s|\le 1/2$ we thus obtain from \eqref{g7} for $k,l \ge k_0$ that
\be{g8}
\left|\log y\ol_j - \log y\ok_j \right| \le \left|g\ol_j - g\ok_j\right|
+ \sum_{i=1}^j \left(4|\Delta\ol_i - \Delta\ok_i|^2 + 2|\Delta\ol_i||\Delta\ol_i - \Delta\ok_i|\right),
\ee
We now use the elementary inequality
$$
4|\Delta\ol_i - \Delta\ok_i|^2 + 2|\Delta\ol_i||\Delta\ol_i - \Delta\ok_i| \le 12 \max_{\ell=1, \dots, j}\left|g\ol_\ell - g\ok_\ell\right| \left(|g\ol_i - g\ol_{i-1}| + |g\ok_i - g\ok_{i-1}|\right)
$$
to conclude that there exists a constant $C>0$ depending only on the difference $g(T) - g(0)$ such that
\be{g9}
\left|\log y\ol_j - \log y\ok_j \right| \le C \max_{i=1, \dots, j}\left|g\ol_i - g\ok_i\right|.
\ee
Note that the sequence $y\ok$ is uniformly bounded
from below by virtue of \eqref{g6}. To derive an upper bound for $k\in\nat$ and $j=1,\dots, m$ we notice that
$$
\log y\ok_j=\sum_{i=1}^j \log\left(1+\frac{\Delta\ok_i}{1- \Delta\ok_i}\right)\le \sum_{i=1}^j \frac{\Delta\ok_i}{1- \Delta\ok_i}\le 2 (g\ok(T)-g\ok(0)) = 2 (g(T)-g(0)).
$$
Hence, $1\le y\ok(t) \le \bar{C}$ for every $t\in[0,T]$ and $k\in\nat$, where
\be{g9a}
\bar{C}= \expe^{2 (g(T)-g(0))}.
\ee
Finally, the uniform convergence follows from \eqref{g9}. Uniqueness is a consequence of the following Gronwall-type statement.
\epf

\begin{lemma}\label{gl2}
Let $g:[0,T] \to \real$ and $y: [0,T] \to [1,\infty)$ be as in Lemma \ref{gl1}, and assume that a right-continuous function $z: [0,T] \to [0,\infty)$ satisfies for some $\gamma \ge 0$ the inequality
\be{g1z}
z(t) \le\gamma + \int_0^t z(\tau)\dd g(\tau) \ \ \forall t\in [0,T].
\ee
Then $z(t) \le \gamma y(t)$ for all $t \in [0,T]$.
\end{lemma}

\bpf{Proof}
Put $v(t) = z(t) - \gamma y(t)$ for $t \in [0,T]$. Then
\be{g1v}
v(t) \le \int_0^t v(\tau)\dd g(\tau) \ \ \forall t\in [0,T]
\ee
and assume that the set $A = \{t\in [0,T]: v(t) > 0\}$ is non-empty. Put $t_0 = \inf A$. Then either $t_0 = 0$ and $v(0) = 0$ by right-continuity, or $t_0 >0$ and,
by Theorem \ref{hake},
$$
0 \le v(t_0) \le \int_0^{t_0} v(\tau)\dd g(\tau) \le v(t_0)(g(t_0)-g(t_0-)),
$$
yielding $\big(1-(g(t_0)-g(t_0-))\big)v(t_0)\le 0$.
Recalling \eqref{g-12} , we conclude that $v(t) = 0$ for all $t \in [0,t_0]$. We now choose a sequence $\{t_n; n\in \nat\}$, $t_n \searrow t_0$ as $n \to \infty$ and such that $v(t_n) > 0$. Let $A_n$ be the sets
$$
A_n = \{t \in [t_0, t_n]: v(t) \ge v(t_n)\},
$$
and put $\hat t_n = \inf A_n$. Then for all $n\in \nat$ we have $v(\hat t_n) \ge v(t_n)$, $\hat{t}_n>t_0$, and, by Lemma \ref{kl9},
$$
v(\hat t_n) \le \int_{t_0}^{\hat t_n} v(\tau)\dd g(\tau) \le v(\hat t_n) (g(\hat t_n) - g(t_0)),
$$
which is a contradiction for $n$ sufficiently large.
\epf

The solution of \eqref{g1} is also known as the generalized exponential function, see \cite{MS}. In particular, for $g$ continuous, the solution is given by $y(t)=\expe^{g(t)-g(0)}$.

\end{document}